\documentclass[twoside,12pt]{article}
\usepackage{amssymb}
\usepackage{amsmath}
\usepackage{amsfonts}
\usepackage[usenames,dvipsnames]{color}

\definecolor{viola}{rgb}{0.3,0,0.7}
\definecolor{ciclamino}{rgb}{0.5,0,0.5}
\definecolor{VioletRed}{rgb}{0.8,0.,0.8}
\newenvironment{bettirev}{\color{blue}}{\color{black}}

\newcommand{\bber}{\begin{bettirev}}
\newcommand{\eber}{\end{bettirev}}
\newcommand{\RR}{\mathbb{R}}

\newcommand{\NN}{\mathbb{N}}
\numberwithin{equation}{section}

\begin{document}

\pagestyle{myheadings}
\newcommand\testopari{\sc \ Colli \ --- \ Marinoschi \ --- \ Rocca}
\newcommand\testodispari{\sc  Sharp interface control in a Penrose-Fife model}
\markboth{\testodispari}{\testopari}

\title{Sharp interface control in a Penrose-Fife model}
\author{Pierluigi Colli \\
Dipartimento di Matematica, Universit\`a degli Studi di Pavia,\\
Via Ferrata~1, 27100 Pavia, Italy\\
E-mail: \texttt{pierluigi.colli@unipv.it} \\
\and Gabriela Marinoschi \\
Institute of Mathematical Statistics and Applied Mathematics,\\
 Calea 13 Septembrie No.13, Sector 5, 050711, Bucharest, Romania\\
E-mail: \texttt{gmarino@acad.ro} \and Elisabetta Rocca \\
Weierstrass Institute for Applied Analysis and Stochastics, \\
Mohrenstr.~39, D-10117 Berlin, Germany\\
E-mail: \texttt{elisabetta.rocca@wias-berlin.de}\\
and\\
Dipartimento di Matematica, Universit\`a degli Studi di Milano,\\
Via Sal\-di\-ni 50, 20133 Milano, Italy\\
E-mail \texttt{elisabetta.rocca@unimi.it}}
\date{}
\maketitle

\vspace{-5mm}

\begin{abstract}
In this paper we study a singular control problem for
a system of PDEs describing a phase-field model of Penrose-Fife type. The
main novelty of this contribution consists in the idea of forcing a sharp
interface separation between the states of the system by using heat sources
distributed in the domain and at the boundary. We approximate
the singular cost functional with a regular one, which is based on the
Legendre-Fenchel relations. Then, we obtain a regularized control problem for
which we compute the first order optimality conditions using an adapted
penalization technique. The proof of some convergence results and the
passage to the limit in these optimality conditions lead to the
characterization of the desired optimal controller.

\vspace{2mm}


\vspace{2mm}

\noindent \textbf{AMS (MOS) subject clas\-si\-fi\-ca\-tion:} 49J20, 82B26,
90C46.
\end{abstract}

~~


\section{Introduction}

We are concerned with a control problem of a system governed by the
Penrose-Fife phase transition model. Using the distributed heat source and
the boundary heat source as controllers we aim at forcing a sharp interface
separation between the states of the system, while keeping its temperature
at a certain
average
level $\theta _{f}$.

The phase-field model considered here has been proposed by Penrose and Fife
in \cite{P-F-1} and \cite{P-F-2} as a thermodynamically consistent model for
the description of the kinetics of phase transition and phase separation
processes in binary materials. It is a PDE system coupling a singular heat
equation (as seen in (\ref{PF1}) below) for the absolute temperature $\theta
$ with a nonlinear equation which describes the evolution of the phase
variable $\varphi $ (see {(\ref{PF2})}), which represents the local
fraction of one of the two components. These equations are accompanied by
initial data for $\theta $ and $\varphi $  (cf. \eqref{PF5} and \eqref{PF6}) and by boundary conditions,
considered here of Robin type for $\theta $  (cf. \eqref{PF3}) and of homogeneous Neumann type
for $\varphi$  (cf. \eqref{PF4}), according to physical considerations. As far as the
Penrose-Fife model is concerned, a vast literature is devoted to the
well-posedness (cf., e.g.~{\cite{ColliLau98, ColliLauSpr, HSZ,
laurencot, Lau, MRSS, Sch-S-Z}}) and to the long-time behavior of solutions
both in term of attractors (cf., e.g.~{\cite{IK, R-Sch,Sc09}}) and of
convergence of single trajectories to stationary states (cf., e.g.,~{\cite{CHIS, 
FS}}), while the associated control problem is less
studied in the literature.

A control problem was introduced first in \cite{SprZh} for a Penrose-Fife
type model with Robin-type boundary conditions for the temperature and a
heat flux proportional to the gradient of the {inverse} absolute
temperature. The first order optimality conditions were derived without
imposing any local constraint on the state and only in case of a double-well
potential in the phase equation. Later on the study has been refined in \cite%
{H-S-S}, where the authors succeeded in removing {such restrictions} on
the problem and treating the case with state constraints.

Let us finally quote the paper \cite{F-F-M} where a phase transition system
was controlled by means of the heat supply in order to be guided into a
certain state with a solid (or liquid) part in a prescribed subset of the
space domain{, by maintaining} it in this state during a period of
time. The system was controlled to form a diffusive boundary between the
solid and liquid states.

Coming back to our problem, we assume here
that the phase transition takes place in the interval $(0,T),$ with $T$
finite{, and that} the system occupies an open bounded domain $\Omega $
of $\mathbb{R}^{3},$ having the boundary $\Gamma $ sufficiently smooth. The
Penrose-Fife system we are interested in reads (see~{\cite{P-F-1, 
ColliLau98, R-Sch}})
\begin{equation}
\theta_{t}-\Delta \beta (\theta )+\varphi _{t}=u,\ \mbox{ in }Q:=(0,T)\times
\Omega ,  \label{PF1}
\end{equation}
\begin{equation}
\varphi _{t}-\Delta \varphi +(\varphi ^{3}-\varphi )=\frac{1}{\theta _{c}}-%
\frac{1}{\theta },\ \mbox{ in }Q,  \label{PF2}
\end{equation}
\begin{equation}
-\frac{\partial \beta (\theta )}{\partial \nu }=\alpha (x)(\beta (\theta
)-v),\ \mbox{ on }\Sigma :=(0,T)\times \Gamma ,  \label{PF3}
\end{equation}
\begin{equation}
\frac{\partial \varphi }{\partial \nu }=0,\ \mbox{ on }\Sigma ,  \label{PF4}
\end{equation}
\begin{equation}
\left. \theta \right\vert _{t=0}=\theta _{0},\ \mbox{ in }\Omega ,  \label{PF5}
\end{equation}
\begin{equation}
\left. \varphi \right\vert _{t=0}=\varphi _{0},\ \mbox{ in }\Omega ,
\label{PF6}
\end{equation}
where $\beta \in C^{1}(0,\infty )$
{and $\beta(r)$ behaves like $- c_1\,/r $
closed to $0$ and like $c_2\, r $ in a neighborhood of $+\infty$, for some constants $c_1$ and $c_2$. Then, for the sake of simplicity we can assume that}
\begin{equation}
\beta (r)={-\frac{1}{r}+ r.}  \label{PF8}
\end{equation}
{Next, we let}
\begin{equation}
{\alpha \in H^1(\Gamma) \cap L^{\infty}(\Gamma),}\quad 0<\alpha
_{m}\leq \alpha (x)\leq \alpha _{M}\ {\mbox{ a.e. }x\in \Gamma},
\label{PF9}
\end{equation}
with $\alpha _{m},$ $\alpha _{M}$ constants. The constant $\theta _{c}$ is the 
transition temperature, $u$ is the distributed heat source 
and $v$ is the boundary heat source.

Note that the heat flux law \eqref{PF8} is a common choice in several types
of phase-transition and phase-separation models both in liquids and in
crystalline solids (cf., e.g., {\cite{ColliLau98, P-F-1, R-Sch}} where
similar growth conditions are {postulated}).

We denote the Heaviside function (translated by $\theta _{c})$ by
\begin{equation}
H(r)=\left\{
\begin{array}{l}
1,\mbox{ \ \ \ \ \ \ \ }r>\theta _{c} \\
\lbrack -1,1],\ r=\theta _{c} \\
-1,\mbox{ \ \ \ \ \ }r<\theta _{c}
\end{array}
\right. \label{defH}
\end{equation}
{and this will be useful to set the third control variable in our
problem. Indeed, let} us define
the cost functional as
\begin{equation}
J(u,v,\eta )=\lambda _{1}\int_{Q}(\theta -\theta_{f})^{2}\,{dx\hskip1pt
dt} +\lambda _{2}\int_{Q}(\varphi -\eta )^{2}\,{dx\hskip1pt dt}
\label{J}
\end{equation}
{and introduce the control problem:}
\begin{equation}
\mbox{Minimize }J(u,v,\eta )\mbox{ for all }(u,v,\eta )\in K_{1}\times
K_{2}\times K_{3},  \tag{$P$}
\end{equation}
subject to (\ref{PF1})--(\ref{PF6}), where
\begin{equation}
K_{1}=\{u\in L^{\infty }(Q)\ {:}\ \ u_{m}\leq u(t,x)\leq u_{M}\
\mbox{
a.e. }(t,x)\in Q\},  \label{K1}
\end{equation}
\begin{equation}
K_{2}=\{v\in L^{\infty }(\Sigma )\ {:}\ \ v_{m}\leq v(t,x)\leq v_{M}\
\mbox{
a.e. }(t,x)\in \Sigma \},  \label{K2}
\end{equation}
\begin{equation}
K_{3}=\{\eta \in L^{\infty }(Q)\ {:}\ \ \eta (t,x)\in H(\theta (t,x))\
\mbox{
a.e. }(t,x)\in Q\},  \label{K3}
\end{equation}
and $u_{m},$ $u_{M},$ $v_{m},$ $v_{M}$ are fixed real values. The
positive constants $\lambda _{1},$ $\lambda _{2}$ are used to give
{more importance to one term or the other} in $(P).$

{With a general approach,} we can consider that
\begin{equation}
{\mbox{$\theta_{f}$\ is a function of $t$\
and $x,$\ and }\, \theta_{f}\in L^{2}(Q) .}  \label{pier1}
\end{equation}
All the results in this paper {hold under this condition}. If by the
control problem one intends to preserve the system separated in two phases
by the sharp interface{,} it should be added that $\theta_{f}$\ must
belong to a neighboorhood of $\theta _{c},$\ i.e., $\left\Vert
\theta_{f}-\theta _{c}\right\Vert _{L^{2}(Q)}\leq \delta ,$\ with $\delta $\ %
{rather} small.

 The problem $(P)$ is introduced in order to enforce the formation of a sharp interface between the two
phases by the constraint $\eta \in K_{3}.$ As far {as} we know such a
control problem has not been previously studied.

Let us note, however, that the well-posedness of an initial-boundary value
Stefan-type problem with phase relaxation or with standard interphase
equilibrium conditions (cf. \eqref{K3}), where the heat flux is proportional
to the gradient of the inverse absolute temperature, was studied in \cite%
{CS95} and \cite{CS97}, for Robin-type boundary conditions. It was shown in
these contributions that the Stefan problems with singular heat flux are the
natural limiting cases of a thermodynamically consistent model of
Penrose-Fife type.

The layout of this paper is as follows. In Section 2, Theorem 2.2, we review
the existence results for the state system and provide new results
concerning the supplementary regularity of the state which will be necessary
in the computation of the optimality conditions. Then, we prove in Theorem
2.3 the existence of at least {one} solution to problem $(P),$
represented by an optimal triplet of controllers $(u,v,\eta )$ and the
corresponding pair of states $(\theta ,\varphi ).$

Due to the singularity induced by the graph representing the sharp
interface, the conditions of optimality cannot be deduced directly for $(P).$
In order to avoid working with the graph $H(\theta )${, in
Section 3 we introduce} an approximating problem $(P_{\varepsilon })$ in
which the constraint $\eta \in H(\theta )$ is replaced by an equivalent
relation based on the Legendre-Fenchel relations between a proper convex
lower semicontinuous (l.c.s.) function $j$ and its conjugate, $j^{\ast }$.
In this case $j$ is the potential of $H.$ This approximating problem has at
least one solution (see Proposition 3.1) which is the appropriate
approximation of a solution to $(P).$ This last assertion relies on the
convergence result of $(P_{\varepsilon })$ to $(P)$ given in Theorem 3.2. In
Section~4, we rigorously examine the question concerning the computation of
the optimality conditions. A second approximation is represented by a
penalized minimization problem $(P_{\varepsilon ,\sigma })$ in which $j$ is
replaced by its {Moreau-Yosida} regularization. The optimality conditions for $
(P_{\varepsilon ,\sigma })$ are provided by explicit expressions in
Proposition~4.5. Some estimates and the proof of the strong convergence (as $
\sigma \rightarrow 0)$ of the controllers in $(P_{\varepsilon ,\sigma })$
allow the {passage} to the limit as $\sigma \rightarrow 0$ in order to
recover the form of the controllers in problem $(P_{\varepsilon })${ : this is 
performed in Theorem~4.6. Recalling~Theorem 2.3, the optimal controller in $(P)$ is
obtained as the limit of a sequence of optimal controllers in $
(P_{\varepsilon })$, on the basis of the convergence of $(P_{\varepsilon })$
to $(P).$

\section{Existence in the state system and control problem}

We denote by $V$ the Sobolev space $H^{1}(\Omega )$ endowed with the
standard scalar product%
\begin{equation}
\left\Vert \psi \right\Vert _{V}=\left( \int_{\Omega }\left\vert \nabla \psi
(x)\right\vert ^{2}dx+\int_{\Omega }\left\vert \psi (x)\right\vert ^{2}
dx\right)^{\!1/2}.  \label{H1-norm}
\end{equation}
{We identify $L^2(\Omega)$ with its dual space, in order that $V
\subset L^2(\Omega) \subset V'$ with dense and compact embeddings.} We
recall that if
$\alpha$ satisfies \eqref{PF9},
then the norm
\begin{equation}
\vert\vert\vert\, \psi \, \vert\vert\vert =\int_{\Omega }\left\vert \nabla
\psi (x)\right\vert ^{2}dx+\int_{\Gamma }\alpha{(x)}\left\vert \psi
(x)\right\vert ^{2}{ds}  \label{H1-necas}
\end{equation}
is equivalent {to} $\left\Vert \psi \right\Vert _{V},$ due to the
inequality
\begin{equation}
\left\Vert \psi \right\Vert _{V}^{2}\leq C_{P}\left( \int_{\Omega
}\left\vert \nabla \psi (x)\right\vert ^{2}dx+\int_{\Gamma }\left\vert \psi
(x)\right\vert ^{2}ds \right) ,\ \forall \psi \in V,  \label{necas}
\end{equation}
(see \cite[{p.~20}]{necas-67}), with $C_{P}$ depending on $\Omega .$
For simplicity, {in the following let us not indicate the arguments of functions} in the integrals.

\medskip

\noindent \textbf{Definition 2.1.} Let
\begin{align}
\theta _{0} \in L^{2}(\Omega ),\quad \theta _{0}>0\ \mbox{ a.e. in } \Omega
,\quad \ln \theta _{0}\in L^{1}(\Omega ),\quad \varphi _{0} \in H^{1}(\Omega
),  \notag \\
u\in L^{2}(Q),\quad v\in L^{2}(\Sigma ), \quad
{\alpha \mbox{ satisfies
\eqref{PF9}}}.  \label{10}
\end{align}
We call a \textit{solution} to (\ref{PF1})--(\ref{PF6}) a pair $(\theta
,\varphi )$ such that
{
\begin{align}\label{11}
&\theta  \in L^{2}(0,T;V) \cap C([0,T];L^{2}(\Omega ))\cap W^{1,2}([0,T];V^{\prime }),
 \quad \beta (\theta ),\, \frac1\theta  \in L^{2}(0,T;V), \\
\label{12}
&\varphi \in L^{2}(0,T;H^{2}(\Omega )) \cap C([0,T];V)\cap
W^{1,2}([0,T];L^{2}(\Omega )),
\end{align}which satisfies (\ref{PF1})--(\ref{PF4}) in the form
\begin{align}\nonumber
\int_{0}^{T}\left\langle \frac{d\theta }{dt}(t),\psi _{1}(t)\right\rangle
_{V^{\prime },V}dt+\int_{Q}\nabla \beta (\theta )\cdot \nabla \psi
_{1}\,{dx\hskip1pt dt} +\int_{Q}\frac{d\varphi }{dt}\psi _{1}\,{dx\hskip1pt dt}  \\
\label{13}
 +\int_{\Sigma }\alpha \beta (\theta )\psi _{1}\,{ds\hskip1pt dt} =
 \int_{Q}u\psi_{1}\,{dx\hskip1pt dt} +\int_{\Sigma }\alpha v\psi _{1}\,{ds\hskip1pt dt} ,
\end{align}
\begin{align}
\nonumber
\int_{Q}\frac{d\varphi }{dt}\psi _{2}\,{dx\hskip1pt dt} +\int_{Q}\nabla \varphi \cdot
\nabla \psi _{2}\,{dx\hskip1pt dt}
+ \int_{Q}(\varphi ^{3}-\varphi )\psi _{2}\,{dx\hskip1pt dt}  \\
 \label{14}
=\int_{Q}\left( \frac{1}{\theta _{c}}-\frac{1}{\theta }\right) \psi
_{2}\,{dx\hskip1pt dt} ,
\end{align}
for any $\psi _{1},$ $\psi _{2}\in L^{2}(0,T;V),$ and such that
the initial conditions (\ref{PF5})--(\ref{PF6}) hold}. \medskip

{The next statement collects a number of properties of solutions to
(\ref{PF1})--(\ref{PF6}).} \medskip

\noindent \textbf{Theorem 2.2. }\textit{Let} \textit{assumptions }(\ref{10})
\textit{hold. Then }(\ref{PF1})--(\ref{PF6}) \textit{has a unique solution, 
{fulfiling}}
\begin{equation*}
{\theta >0\ \mbox{ \textit{{a.e. in}} }Q, \quad \ln \theta \in
L^{\infty }(0,T;L^{1}(\Omega ))}
\end{equation*}
\textit{{and} satisfying the estimates}
\begin{equation}
{\left\Vert \theta \right\Vert _{L^{2}(0,T;V)} + \left\Vert \theta
\right\Vert _{L^{\infty }(0,T;L^{2}(\Omega ))}} +\left\Vert \theta
\right\Vert _{W^{1,2}([0,T];V^{\prime })}  +\left\Vert\frac{1}{\theta}\right\Vert_{L^{2}(0,T;V)}\leq C, \label{15}
\end{equation}
\begin{equation}
{\left\Vert \varphi \right\Vert _{L^{2}(0,T;H^{2}(\Omega ))}
+\left\Vert \varphi \right\Vert _{L^{\infty }(0,T;V)}} +\left\Vert \varphi
\right\Vert_{W^{1,2}([0,T];L^{2}(\Omega ))}\leq C.  \label{16}
\end{equation}
\textit{Moreover, let us set $\bar\theta:=\theta_1-\theta_1$, $
\bar\varphi:=\varphi_1-\varphi_2$, $\bar u:=u_1-u_2$, $\bar v:=v_1-v_2$,
where $(\theta_1, \varphi_1)$, and $(\theta_2, \varphi_2)$ are the solutions
of \eqref{PF1}--\eqref{PF6} corresponding respectively to the data $u_1,\,
v_1$ and $u_2,\,v_2$, to the same initial data $\theta_0$, $\varphi_0$ and
to the same coefficient $\alpha$; then, we have the following continuous
dependence estimate of the solution with respect to the data:}
\begin{align}  \label{cd}
\|\bar\theta\|_{L^2(Q)}^2+\|\bar\varphi\|_{C([0,T];L^2(\Omega))}^2+\|\bar%
\varphi\|_{L^2(0,T;V)}^2\leq C\left(\|\bar u\|_{L^2(Q)}+\|\bar
v\|_{L^2(\Sigma)}^2\right) ,
\end{align}
\textit{with the positive constant $C$ depends only on the problem %
{parameters}, but not on $u_i, \, v_i$, {$1=1,2$.}}
\textit{{Next, we list some regularity properties
of the solution: if,} in addition to \eqref{10}, {we suppose that}}
\begin{equation}
\varphi _{0}\in H^{2}(\Omega ),\ \ \frac{\partial \varphi _{0}}{\partial \nu
}=0\mbox{ on }\Gamma ,  \label{17}
\end{equation}
\textit{then {we have}}
\begin{equation}
\varphi \in { L^{\infty }(Q)\cap L^{\infty }(0,T;H^{2}(\Omega )) \cap
W^{1,2}([0,T];V)}  \label{18}
\end{equation}
\textit{and }
\begin{equation}
\left\Vert \varphi \right\Vert _{L^{\infty }(Q)}+\left\Vert \varphi
\right\Vert _{L^{\infty }(0,T;H^{2}(\Omega ))}+\left\Vert \varphi
\right\Vert _{W^{1,2}([0,T];{V})}\leq C;  \label{19}
\end{equation}

\smallskip\noindent \textit{{further,} if, in addition to \eqref{10}, %
{there hold}}
\begin{align}
{\theta _{0} , \, \frac{1}{\theta _{0}}\in L^{\infty }(\Omega ), \quad
u \in L^{2}(0,T;L^{6}(\Omega )), \quad v\in L^{\infty }(\Sigma )\
\,\hbox{\it and}\qquad}  \notag \\
{ v\leq v_M \quad \mbox{\it a.e. in} \ \, \Sigma, \label{20}}
\end{align}
{\textit{then we have }\begin{equation}
\theta , \,   \frac{1}{\theta }\in L^{\infty }(Q)
\label{21}
\end{equation}\textit{with }\begin{equation}
\left\Vert \theta \right\Vert _{L^{\infty }(Q)}+\left\Vert \frac{1}{\theta }\right\Vert _{L^{\infty }(Q)}\leq C,  \label{21-1}
\end{equation}}\textit{where }$C$\textit{\ denotes several positive
constants depending {only} on the problem parameters.}

\medskip

\noindent \textbf{Proof. }
The proof of existence of solutions to (\ref{PF1})--(\ref{PF6}) follows from
an adaptation of \cite[Thm.~2.3]{ColliLau98} to the case {of} $\alpha$
non constant in \eqref{PF3}. The uniqueness of solutions has been proved in
\cite[Thm.~1]{ColliLauSpr} and it has been then generalized to the case of
less regular data (satisfying exactly assumptions \eqref{10}) in \cite[%
Thm.~3.5]{R-Sch}, where also a continuous dependence result of the solution
with respect to the data has been
{shown. We also refer to the
above-mentioned papers for the proof of estimates \eqref{15}--\eqref{16}.} In
what follows the positive constants $C$, which may also differ from line to
line, will depend only on the problem parameters.

\paragraph{Proof of estimate \eqref{cd}. } Following the lines of \cite[Thm.~1]{ColliLauSpr} and \cite[Thm.~3.5]{R-Sch},
we write \eqref{PF1} firstly for $(\theta_1, \varphi_1)$ and then for $
(\theta_2, \varphi_2)$, being $(\theta_1, \varphi_1)$, and $(\theta_2,
\varphi_2)$ two solutions (\ref{PF1})--(\ref{PF6}) corresponding
respectively to the data $u_1,\, v_1$ and $u_2,\,v_2$, to the same initial
data $\theta_0$, $\varphi_0$, and to the same coefficient $\alpha$. Taking
the difference {and integrating with respect to time, we easily obtain}
\begin{align}  
&\int_{0}^{t}\left\langle \bar \theta (\tau),\psi _{1}(\tau)\right\rangle
_{V^{\prime },V}d\tau+\int_{Q_t}1*\nabla\left( \beta (\theta_1
)-\beta(\theta_2)\right)\cdot \nabla \psi _{1}\,{dx\hskip1pt d\tau} \notag \\
&\qquad +\int_{Q_t}\bar \varphi \psi _{1}\,{dx\hskip1pt d\tau}  
+\int_{\Sigma_t }\alpha \left(1*\left(\beta (\theta_1
)-\beta(\theta_2)\right)\right)\psi _{1}\,{ds\hskip1pt d\tau} \notag \\
&{}=
\int_{Q_t}\left(1*\bar u\right)\psi_{1}\,{dx\hskip1pt d\tau} +\int_{\Sigma_t
}\alpha \left(1*\bar v\right)\psi _{1}\,{ds\hskip1pt d\tau} ,
\label{13int}
\end{align}
for any $\psi _{1}\in L^{2}(0,T;V),$ where $Q_t:=(0,t) \times \Omega $, $
\Sigma_t:= (0,t) \times \Gamma $. Here, we have also used the notation $
\bar\theta:=\theta_1-\theta_1$, $\bar\varphi:=\varphi_1-\varphi_2$, $\bar
u:=u_1-u_2$, $\bar v:=v_1-v_2$, and denoted by $*$ the standard time
convolution operator,
{so that
$\displaystyle (a* b) (t) = \int_0^t a(t-\tau)b(\tau)d\tau$, $t\in (0,T]$.}
Choosing now as test function $\psi_1=\beta(\theta_1)-\beta(\theta_2)$, and
using the monotonicity properties of the function $\theta\mapsto-1/\theta$,
we find out that
\begin{align}
&\int_{Q_t}|\bar\theta|^2\,{dx\hskip1pt d\tau}+\frac12\int_\Omega|\nabla
1*\left(\beta(\theta_1)-\beta(\theta_2)\right)(t)|^2\, dx  \notag \\
&\qquad +\int_{Q_t}\bar\varphi\left(\beta(\theta_1)-\beta(\theta_2)\right)\,
{dx\hskip1pt d\tau}+\frac12\int_\Gamma
\alpha\left|1*\left(\beta(\theta_1)-\beta(\theta_2)\right)(t)\right|^2\, ds
\notag \\
&\leq \int_{Q_t}(1*\bar u)\left(\beta(\theta_1)-\beta(\theta_2)\right) \,{dx%
\hskip1pt d\tau} +\int_{\Sigma_t}\alpha \left(1*\bar
v\right)\left(\beta(\theta_1)-\beta(\theta_2)\right)\,{ds\hskip1pt d\tau}.
\label{cd1}
\end{align}
{Next}, taking the differences of \eqref{PF2}, testing by $\psi_2=\bar\chi$,
and exploiting the monotonicity of $\varphi\mapsto \varphi^3$, we have that
\begin{align}  \label{cd2}
&\frac12\int_\Omega|\bar\varphi(t)|^2\,dx+\int_{Q_t} |\nabla\bar\varphi|^2\,
{dx\hskip1pt d\tau}\leq \int_{Q_t}|\bar\varphi|^2\,{dx\hskip1pt d\tau}
+\int_{Q_t}\left(-\frac{1}{\theta_1}+\frac{1}{\theta_2}\right)\bar\varphi\, {%
dx\hskip1pt d\tau}.
\end{align}
Now, summing up \eqref{cd1} and \eqref{cd2}, we take advantage of a
cancellation of one term due to the special form \eqref{PF8} of $\beta$.
Then, in view of assumption \eqref{PF9} on $\alpha$, we arrive~at
\begin{align}
&\int_{Q_t}|\theta|^2\,{dx\hskip1pt d\tau}+\|1*\left(\beta(\theta_1)-\beta(%
\theta_2)\right)(t)\|_V^2+\|\bar\varphi(t)\|_{L^2(\Omega)}^2+\int_{Q_t}|%
\nabla\bar\varphi|^2\,{dx\hskip1pt d\tau}  \notag \\
&\leq C_1\left(\int_{Q_t}|\bar\varphi\bar\theta|\,{dx\hskip1pt d\tau}
+\int_{Q_t}|(1*\bar u)\left(\beta(\theta_1)-\beta(\theta_2)\right) |\,{dx%
\hskip1pt d\tau}\right.  \notag \\
&\qquad\quad\left.+\int_{\Sigma_t}|\alpha \left(1*\bar
v\right)\left(\beta(\theta_1)-\beta(\theta_2)\right)|\,{ds\hskip1pt d\tau}
+\int_{Q_t}|\bar\varphi|^2\,{dx\hskip1pt d\tau}\right).  \label{cd3}
\end{align}
Let us now estimate the integrals on the right hand side of \eqref{cd3} as
follows. Using Young's inequality, we deduce that
\begin{align}  \label{cd4}
\int_{Q_t}|\bar\varphi\bar\theta|\,{dx\hskip1pt d\tau}&\leq
\delta_1\int_{Q_t}|\bar \theta|^2\,{dx\hskip1pt d\tau}+C_{\delta_1}
\int_{Q_t}|\bar\varphi|^2\,{dx\hskip1pt d\tau},
\end{align}
for some $\delta_1>0$ to be chosen later. Finally,  integrating by parts in
time and using again Young's inequality, we obtain
\begin{align}
&\int_{Q_t}|(1*\bar u)\left(\beta(\theta_1)-\beta(\theta_2)\right) |\, {dx%
\hskip1pt d\tau}  \notag \\
&\leq \int_\Omega|1*\bar u|(t)|1*\left(\beta(\theta_1)-\beta(\theta_2)
\right) |(t)\,dx + \int_{Q_t}|\bar
u\left(1*\left(\beta(\theta_1)-\beta(\theta_2)\right)\right) |\,{dx\hskip1pt
d\tau}  \notag \\
&\leq \delta_2\|1*\left(\beta(\theta_1)-\beta(\theta_2)\right) (t)\|_{V}^2
+C_{\delta_2}\|1*\bar u(t)\|_{L^2(\Omega)}^2  \notag \\
&\qquad{}+ \int_{Q_t}|\bar u|^2\,{dx\hskip1pt d\tau}+
\int_0^t\|1*\left(\beta(\theta_1)-\beta(\theta_2)\right)\|_{L^2(\Omega)}^2%
\,d\tau  \label{cd5}
\end{align}
and, thanks to \eqref{PF9},
\begin{align}
& \int_{\Sigma_t}|\alpha \left(1*\bar
v\right)\left(\beta(\theta_1)-\beta(\theta_2)\right)|\,{ds\hskip1pt d\tau}
\notag \\
&\leq \int_\Gamma|\alpha||1*\bar
v|(t)|1*\left(\beta(\theta_1)-\beta(\theta_2)\right) |(t) \,ds +
\int_{\Sigma_t} |\alpha \bar
v\left(1*\left(\beta(\theta_1)-\beta(\theta_2)\right)\right)|\,{ds\hskip1pt
d\tau}  \notag \\
&\leq \delta_3\|1*\left(\beta(\theta_1)-\beta(\theta_2)\right) (t)\|_{V}^2
+C_{\delta_3}\|1*\bar v(t)\|_{L^2(\Gamma)}^2  \notag \\
&\qquad {}+ \int_{\Sigma_t}|\bar v|^2\,{ds\hskip1pt d\tau}
+\int_0^t\|1*\left(\beta(\theta_1)-\beta(\theta_2)\right)\|_{L^2(\Gamma)}^2\,
{d\tau} .  \label{cd6}
\end{align}
Collecting now estimates (\ref{cd3})--(\ref{cd6}) and choosing the constants
$\delta_i$, $i=1,2,3,$ such that $C_1(\delta_1+\delta_2+\delta_3)<1$, we
infer that
\begin{align*}
&\int_{Q_t}|\theta|^2\,{dx\hskip1pt d\tau}+\|1*\left(\beta(\theta_1)-\beta(%
\theta_2)\right)(t)\|_V^2+\|\bar\varphi(t)\|_{L^2(\Omega)}^2+\int_{Q_t}|%
\nabla\bar\varphi|^2\,{dx\hskip1pt d\tau}  \notag \\
&\leq C_2\left(\int_0^t\|\bar\varphi(\tau)\|_{L^2(\Omega)}^2
d\tau+\int_0^t\|1*\left(\beta(\theta_1)-\beta(\theta_2)\right)\|_{V}^2\,{%
d\tau} + \|\bar u\|^2_{L^2(Q_t)} + \|\bar v\|^2_{L^2(\Sigma_t)} \right),
\end{align*}
from which, using a standard Gronwall lemma, we deduce the desired \eqref{cd}.

\paragraph{Proof of estimate \eqref{19}.} 
In order to prove the regularity \eqref{18} and estimate \eqref{19}, we can
proceed formally
{testing \eqref{PF2} by  $-\Delta \varphi_t$ and
integrating by parts with the help of \eqref{PF4}.} This choice should be
made rigorous in the framework of a regularized scheme, e.g., of
Faedo-Galerkin type, but we prefer to proceed formally here in order not to
overburden the presentation. Making this {formal computation} and
integrating the resulting equation over $(0,t)$, $t\in (0,T]$, we get
\begin{align}
&\frac12\|\Delta
\varphi(t)\|_{L^2(\Omega)}^2-\frac12\|\Delta\varphi_0\|_{L^2(\Omega)}^2+%
\int_0^t\|\nabla\varphi_t\|_{L^2(\Omega)}^2  \, d\tau\notag \\
&\leq \int_0^t\left\|\frac{1}{
\theta}\right\|_V\|\nabla\varphi_t\|_{L^2(\Omega)}  \, d\tau+\int_0^t\int_\Omega\nabla (\varphi^3-\varphi)\nabla \varphi_t\,{dx\hskip1pt d\tau}\,.  \label{e1}
\end{align}
In order to estimate the first integral on the right hand side we can just
use {the} Young inequality together with estimate \eqref{15}. The last
integral, instead, can be treated as follows:
\begin{align}  \label{e2}
\int_0^t\int_\Omega\nabla (\varphi^3-\varphi)\nabla \varphi_t\,{dx\hskip1pt d\tau} \leq C\int_0^t
\left(\|\varphi\|_{L^6(\Omega)}^2+1\right)\|\nabla\varphi\|_{L^6(\Omega)}\|%
\nabla\varphi_t\|_{L^2(\Omega)}  \, d\tau\notag \\
\leq {\frac14} \int_0^t\|\nabla\varphi_t\|_{L^2(\Omega)}^2 \, d\tau+{C}
\int_0^t (1+ \|\Delta\varphi\|_{L^2(\Omega)}^2) \, d\tau\,,
\end{align}
where
{the H\"older and Young inequalities have been used together with
the Gagliardo-Nirenberg} inequality (\cite[p.~125]{nir}) and the previous
{estimate~\eqref{16}. By rearranging in \eqref{e1} and using once more \eqref{16} for the boundedness of $\varphi$ in $L^2(0,T;H^2(\Omega))$,
we obtain the estimate
\begin{equation}\nonumber
\|\Delta\varphi\|_{L^\infty(0,T;L^2(\Omega))}+\|\nabla\varphi_t\|_{L^2(Q)}\leq C\,,
\end{equation}
which, together with \eqref{15}--\eqref{16}, the standard elliptic regularity results}
and the continuous embedding of $L^\infty(0,T; H^2(\Omega))$ into $
L^\infty(Q)$ in 3D, gives the desired \eqref{19}.

\paragraph{Proof of estimate \eqref{21-1}.} We aim first 
to prove the $L^\infty(Q)$-bound for $\theta$. In order to do
that{, we use a Moser-type technique. The procedure consists in testing
\eqref{PF1} by $\ (p+1)\theta^p$, \ $p\in (1,\infty)$.} This estimate is
formal {(cf.~\eqref{13});} indeed, in order to perform it rigorously we
would need to introduce a regularized (truncated) system and then pass to
the limit. However, since the procedure is quite standard, we prefer to
perform only the formal estimate here.
{Testing \eqref{PF1} by $\ (p+1)\theta^p$ \ leads to}
\begin{align}
&\frac{d}{dt}\int_\Omega\theta^{p+1}\, dx+\frac{4p}{p+1}\int_\Omega\left|
\nabla\theta^{\frac{p+1}{2}}\right|^2\, dx +\frac{4p(p+1)}{(p-1)^2}
\int_\Omega\left|\nabla\theta^{\frac{p-1}{2}}\right|^2\, dx
+(p+1)\int_{\Gamma}\alpha\theta^{p+1}\, ds  \notag \\
&=(p+1)\int_\Gamma\alpha \theta^{p-1}\, ds +(p+1)\int_\Gamma\alpha
v\theta^p\, ds +(p+1)\int_\Omega {(u- \varphi_t)}\theta^p\, dx\,.
\notag
\end{align}
Using now assumptions \eqref{10} and \eqref{20}, we get
\begin{align}  
&\frac{d}{dt}\int_\Omega\theta^{p+1}\, dx+\frac{4p}{p+1}\int_\Omega\left|%
\nabla\theta^{\frac{p+1}{2}}\right|^2\, dx
+(p+1)\alpha_m\int_{\Gamma}\left|\theta^{\frac{p+1}{2}}\right|^2\, ds 
\notag \\
&\leq(p+1)\alpha_M\int_\Gamma\theta^{p-1}\,
ds+(p+1)\alpha_Mv_M\int_\Gamma\theta^p\, ds+(p+1)\int_\Omega {(u-
\varphi_t)}\theta^p\, dx\,. \label{e3} 
\end{align}
{Owing to} the Young inequality in the form $a\cdot b\leq \epsilon
\frac{{a}^q}{q}+\frac{1}{\epsilon^{q^{\prime }/q}}\cdot \frac{b^{q^{\prime }}
}{q^{\prime }}$, {$
\frac1q + \frac1{q'} = 1 $,} we estimate the two boundary terms as follows:
\begin{align}  \label{e4}
&\alpha_M\int_\Gamma\theta^{p-1}\, ds\leq \frac{\alpha_m}{4}\frac{(p-1)}{p+1}
\int_\Gamma\theta^{p+1}\, ds+\left(\frac{4}{\alpha_m}\right)^{\!\!\frac{(p-1)%
}{2}}\frac{\alpha_M^{\frac{p+1}{2}}}{(p+1)} 2|\Gamma|\,, \\
\label{e5}
&\alpha_M v_M\int_\Gamma\theta^p\, ds \leq \frac{\alpha_m}{4}\frac{p}{(p+1)}
\int_\Gamma\theta^{p+1}\, ds+\left(\frac{4}{\alpha_m}\right)^{\!p}\frac{%
(\alpha_M v_M)^{p+1}}{(p+1)}|\Gamma|\,.
\end{align}
{Thanks to estimates \eqref{e4}--\eqref{e5}}, \eqref{e3} becomes
\begin{align*}
\frac{d}{dt}\int_\Omega\theta^{p+1}\, dx+\delta
\left(\int_\Omega\left|\nabla\theta^{\frac{p+1}{2}}\right|^2\,
dx+\int_\Gamma\left|\theta^{\frac{p+1}{2}}\right|^2\, ds\right) \\
\leq C R^{p+1}+(p+1)\int_\Omega |{u- \varphi_t}|\theta^p\, dx,
\end{align*}
where $\delta:=\min\{\frac32\alpha_m, \frac83\}>0$ {and $C, \, R $} are
independent of $p$. Using then the continuous embedding of ${V={}}
H^1(\Omega)$ into $L^6(\Omega)$ in 3D, we obtain
\begin{equation}  \label{e6}
\frac{d}{dt}\int_\Omega\theta^{p+1}\, dx+\delta^{\prime
}\left(\int_\Omega\theta^{3(p+1)}\, dx\right)^{1/3}\leq C
R^{p+1}+(p+1)\int_\Omega |{u- \varphi_t}|\theta^p\, dx,
\end{equation}
for some $\delta^{\prime }>0$ always independent of $p$. Then,
{as
\eqref{19} entails the boundedness of $\varphi_t$ in $L^2(0,T;L^6(\Omega))$,}
we estimate the last integral as follows:
\begin{align}  \label{e7}
&(p+1)\int_\Omega |{u- \varphi_t}|\theta^p\, dx  \notag \\
&\leq (p+1)\|{u- \varphi_t}\|_{L^6(\Omega)}\left(\int_\Omega%
\theta^{3(p+1)}\, dx\right)^{1/6}\left(\int_\Omega\theta^{\frac32\frac{(p-1)%
}{2}}\, dx \right)^{2/3}  \notag \\
&\leq \frac{\delta^{\prime }}{2}\left(\int_\Omega\theta^{3(p+1)}\,
dx\right)^{1/3}+C_{\delta^{\prime }}(p+1)^2\|{u- \varphi_t}
\|_{L^6(\Omega)}^2\left(\int_\Omega\theta^{\frac34(p-1)} \, d x\right)^{4/3}  \notag
\\
&\leq \frac{\delta^{\prime }}{2}\left(\int_\Omega\theta^{3(p+1)}\,
dx\right)^{1/3}+C_{\delta^{\prime }}(p+1)^2\|{u- \varphi_t}
\|_{L^6(\Omega)}^2\left(\int_\Omega\theta^{p-1} \, d x \right)\,,
\end{align}
where we have also used the inequality $\|\theta^{\frac{p-1}{2}
}\|_{L^{3/2}(\Omega)}^2\leq C(\Omega)\|\theta^{\frac{p-1}{2}
}\|_{L^2(\Omega)}^2$. Choosing now $p=3$ in {\eqref{e6}--\eqref{e7}}
and
{integrating with respect to time, with the help of \eqref{15} and \eqref{20} we obtain
\begin{equation}\label{e8}
\| \theta\|^4_{L^\infty(0,T;L^4(\Omega))}\leq C\left(1+ \int_0^T\|u- \varphi_t \|_{L^6(\Omega)}^2\, d\tau \right).
\end{equation}
In general}, integrating \eqref{e6} from $0$ to $t$, $t\in (0,T]$, and using %
\eqref{e7} and \eqref{20}, we {infer that}
\begin{equation}  \label{e9}
\int_\Omega\theta^{p+1}(t)\, dx\leq
C\left(R^{p+1}+(p+1)^2\sup_{[0,T]}\left(\int_\Omega\theta^{{\frac34}
(p-1)}\, dx \right)^{\!4/3}\right),
\end{equation}
where $C$ depends on the data but not on $p$. {At this point,} we can
introduce the sequence
\begin{align}
&p_0=4  \notag \\
&p_{k+1}=\frac43p_k+2, \quad k\in \mathbb{N}\,  \notag
\end{align}
and {take} $p=p_{k+1}-1$ in \eqref{e9}, getting
\begin{equation*}
\int_\Omega\theta^{p_{k+1}}(t)\, dx\leq
C\left(R^{p_{k+1}}+(p_{k+1})^2\sup_{[0,T]}\left(\int_\Omega\theta^{p_k}\, dx
\right)^{4/3}\right){.}
\end{equation*}
{We can apply} now \cite[Lemma A.1]{Lau93} with the choices $a=4/3$, %
{$b=c=2$,} $\delta_0=4$, $\delta_k=p_k${. Thus, we deduce that}
\begin{equation*}
\sup_{[0,T]}\|\theta\|_{L^{p_k}(\Omega)}\leq C\,,
\end{equation*}
where $C$ is independent of $k$. Hence, letting $k$ tend to $\infty$, we get
\begin{equation}  \label{e10}
\|\theta\|_{L^\infty(Q)}\leq C\,.
\end{equation}

Finally, we aim to prove the $L^\infty(Q)$-bound for $1/\theta$. Hence, let
us call $h=1/\theta$ and rewrite formally {\eqref{PF2}--\eqref{PF3}} as
follows
\begin{align}  \label{ph1}
&h_t-h^2\Delta\left(h-\frac1h\right)=-h^2 {(u-\varphi_t)}{,} \quad %
\hbox{ {in} }Q, \\  \label{ph2}
&-{\frac{\partial}{\partial\nu}}\left(h-\frac1h\right)=\alpha\left(h-%
\frac1h+v\right){,}\quad\hbox{on }\Sigma\,.
\end{align}
Note that, due to the estimate \eqref{e10}, {there exists} a positive
constant $\bar C$ (depending on the data) such that
\begin{equation*}
h {(t,x)}\geq \bar C\quad \hbox{a.e. } {(t,x)}\in \bar Q\,.
\end{equation*}
Test now \eqref{ph1} by $ph^{p-1}$, $p\in (1,\infty)$, getting
\begin{align}
&\frac{d}{dt}\int_\Omega h^p\, dx+\frac{4p(p+1)}{(p+2)^2}\int_\Omega\left|%
\nabla h^{\frac{p+2}{2}}\right|^2\, dx+\frac{4(p+1)}{p}\int_\Omega\left|%
\nabla h^{\frac{p}{2}}\right|^2\, dx+p\int_\Gamma\alpha h^{p+2}\, ds  \notag
\\
&=p\int_\Gamma\alpha h^p\, ds+p\int_\Gamma\alpha v h^{p+1}\, ds-p\int_\Omega %
{(u- \varphi_t)} h^{p+1}\, dx\,.  \notag
\end{align}
Using now the Young inequality and the assumption \eqref{20}, we end up with
\begin{equation*}
\frac{d}{dt}\int_\Omega h^p\, dx+\delta\left(\int_\Omega\left|\nabla h^{%
\frac{p+2}{2}}\right|^2\, dx +\int_\Gamma \left|h^{\frac{p+2}{2}}\right|^2\,
ds\right)\leq CR^{p+2}+p\int_\Omega |{u- \varphi_t}|h^{p+1}\, dx\, {,}
\end{equation*}
where $\delta$ and $R$ are positive constants independent of $p$.
{By
recalling} the continuous embedding of $H^1(\Omega)$ into $L^6(\Omega)$ in
3D along with H\"older's inequality, we {have that}
\begin{align}
&\frac{d}{dt}\int_\Omega h^p\, dx+\delta^{\prime
}\left(\int_\Omega|h|^{3(p+2)}\, dx\right)^{1/3}  \notag \\
&\leq CR^{p+2}+p\int_\Omega |{u- \varphi_t}|h^{p+1}\, dx  \notag \\
&\leq CR^{p+2}+p \|{u- \varphi_t}\|_{L^6(\Omega)}\|h^{\frac{p+2}{2}
}\|_{L^6(\Omega)}\|h^{\frac{p}{2}}\|_{L^{3/2}(\Omega)}  \notag \\
&\leq CR^{p+2}+\frac{\delta^{\prime }}{2}\|h^{\frac{p+2}{2}}\|_{L^6(\Omega)
}^2+\frac{p^2}{2\delta^{\prime }}\|{u- \varphi_t}\|_{L^6(\Omega)}^2\|h^{%
\frac{p}{2}}\|_{L^{3/2}(\Omega)}^2\, .  \notag
\end{align}
Now, integrating over $(0,t)$, $t\in (0,T]$, and using the continuous
embedding of $L^p(\Omega)$ into $L^{3p/4}(\Omega)$ {as well as}
assumption \eqref{20}, we {infer that}
\begin{align}
\int_\Omega h^p(t)\, dx+\delta^{\prime }\left(\int_\Omega|h|^{3(p+2)}\,
dx\right)^{1/3}&\leq C\left(R^{p+2}+p^2\int_0^t\|{u- \varphi_t}
\|_{L^6(\Omega)}^2\|h\|_{L^{\frac{3p}{4}}(\Omega)}^p\, d\tau\right)  \notag
\\
& \leq C\left(R^{p+2}+p^2\int_0^t\|{u- \varphi_t}\|_{L^6(\Omega)}^2\|h%
\|_{L^p(\Omega)}^p\, d\tau\right)\,,  \notag
\end{align}
where $C$ depends on the data, but not on $p$. Choosing now $p=6$ and
applying {the Gronwall lemma}, we obtain the {starting} point for
an iterating procedure which is completely analogous to the one in \cite[%
p.~269]{laurencot}. Hence, we obtain
\begin{equation}  \label{e11}
{\|h\|_{L^\infty(Q)}= {}}\|1/\theta\|_{L^\infty(Q)}\leq C
\end{equation}
and this concludes the proof of Theorem~2.2.
\hfill $\square $

\medskip

The next result proves the existence of a solution to problem $(P).$

\medskip

\noindent \textbf{Theorem 2.3. }\textit{Assume that}\textbf{\ }
\begin{equation}
\theta _{0}\in L^{2}(\Omega ),\quad \theta _{0}>0\ \mbox{ a.e. in }\Omega
,\quad \ln \theta _{0}\in L^{1}(\Omega ),\quad \varphi _{0}\in H^{1}(\Omega
),  \label{10-0}
\end{equation}
\textit{and} (\ref{PF9}) \textit{hold. }\textit{Then }$(P)$\textit{%
\ has at least one solution.}

\medskip

\noindent \textbf{Proof. }Since $J(u,v,\eta )\geq 0,$ it follows that $J$
has an infimum $d$ and {this infimum} is nonnegative. Let $
(u_{n},v_{n},\eta _{n})_{n\geq 1}$ be a minimizing sequence for $J.$\ This
means that $u_{n}\in K_{1},$ $v_{n}\in K_{2},$ $\eta _{n}\in K_{3}$, $
(\theta _{n},\varphi _{n})$ is the solution to (\ref{PF1})--(\ref{PF6})
corresponding to $u_{n},$ $v_{n},$ $\eta _{n},$ and the following
inequalities take place%
\begin{equation}
d\leq \lambda _{1}\int_{Q}(\theta _{n}-\theta _{f})^{2}\,{dx\hskip1pt
dt}+\lambda _{2}\int_{Q}(\varphi _{n}-\eta _{n})^{2}\,{dx\hskip1pt dt}
\leq d+\frac{1}{n},\ n\geq 1.  \label{22}
\end{equation}

Therefore, {possibly taking} subsequences (denoted still by the
subscript {$n$}), we deduce that
{
\begin{align*}
u_{n}\rightarrow u\ \mbox{ {weakly*} in }L^{\infty }(Q),\quad v_{n}\rightarrow v\ \mbox{ {weakly*} in }L^{\infty }(\Sigma ),\qquad
\\
\eta _{n}\rightarrow \eta\  \mbox{ {weakly*} in }L^{\infty }(Q),\ \mbox{ as }
n\rightarrow \infty ,
\end{align*}}and $u\in K_{1},$ $v\in K_{2},$ $\eta \in K_{3}.$
By (\ref{15})--(\ref{16}) we have
\begin{equation*}
\theta _{n}\rightarrow \theta \ \mbox{ weakly in }L^{2}(0,T;V)\cap
W^{1,2}([0,T];V^{\prime }),\mbox{ as }n\rightarrow \infty ,
\end{equation*}
\begin{equation*}
\frac{1}{\theta _{n}}\rightarrow l\ \mbox{ weakly in }L^{2}(0,T;V),%
\mbox{
as }n\rightarrow \infty ,
\end{equation*}
\begin{align*}
\varphi _{n} \rightarrow \varphi \mbox{ weakly in }L^{2}(0,T;H^{2}(\Omega
))\cap W^{1,2}([0,T];L^{2}(\Omega ))\qquad  \\
\mbox{and {weakly*} in }L^{\infty }(0,T;V),\mbox{ as }n& \rightarrow
\infty .
\end{align*}
These {facts imply, by the Aubin-Lions theorem (see
\cite[p.~57]{lions}),} that%
\begin{align*}
& \theta _{n}\rightarrow \theta \ \mbox{ strongly in }L^{2}(0,T;L^{2}(\Omega
)),\mbox{ as }n\rightarrow \infty , \\
& \varphi _{n}\rightarrow \varphi \ \mbox{ strongly in }L^{2}(0,T;V),%
\mbox{
as }n\rightarrow \infty .
\end{align*}
Therefore, on a subsequence {it results that}
\begin{equation*}
\theta _{n}\rightarrow \theta \ \mbox{ a.e. in }Q,\mbox{ as }n\rightarrow
\infty ,
\end{equation*}
{whence}
\begin{equation*}
\frac{1}{\theta _{n}}\rightarrow \frac{1}{\theta }\ \mbox{ a.e. in }Q,%
\mbox{
as }n\rightarrow \infty ,
\end{equation*}
entailing that $l=1/\theta $ a.e. in $Q$ (%
cf., e.g., \cite[Lemme~1.3, p.~12]{lions}).
With the help of the Egorov theorem, we can also conclude that
\begin{equation*}
\frac{1}{\theta _{n}}\rightarrow \frac{1}{\theta }\ \mbox{ strongly in }
{L^{p}(Q),\mbox{ for all $1\leq p <2$,
as }} n\rightarrow \infty .
\end{equation*}
{On the other hand, in view of \eqref{PF8} the above convergences yield}
\begin{equation*}
{\beta(\theta _{n})\rightarrow  \beta (\theta)  \ \mbox{ weakly in }L^{2}(0,T;V) \, \mbox{ and a.e. in } Q,\mbox{ as }n\rightarrow \infty .}
\end{equation*}

Next, since ${\{\varphi _{n}\}}$ is bounded in $L^{\infty }(0,T;V)$ we
deduce that ${\{\varphi _{n}^{3}\}}$ is bounded in $L^{\infty
}(0,T;L^{2}(\Omega ))$ and {consequently}
\begin{equation*}
\varphi _{n}^{3}\rightarrow l_{1}\ \mbox{ {weakly*} in }L^{\infty
}(0,T;L^{2}(\Omega )),\mbox{ as }n\rightarrow \infty .
\end{equation*}
But, there exists a subsequence such that $\varphi _{n}\rightarrow \varphi $ %
{a.e. in} $Q$. This implies that $\varphi _{n}^{3}\rightarrow \varphi
^{3}$ {a.e. in} $Q$ and we conclude that $l_{1}=\varphi ^{3}$ %
{a.e. in} $Q.$

Now, we recall that $\eta _{n}\in H(\theta _{n})$ {}{a.e. in}{} $Q$, $
\theta _{n}\rightarrow \theta $ strongly in $L^{2}(Q)$ and $\eta
_{n}\rightarrow \eta $ {weakly*} in $L^{\infty }(Q).$ On the basis of %
{the} maximal monotonicity of $H${,} we {deduce} that $\eta
\in H(\theta )$ {a.e. in} $Q.$

Moreover, since the trace operator is continuous from $V$ to $L^{2}(\Gamma ),
$ we
{have that
$\left\Vert \beta (\theta _{n})\right\Vert _{L^2(0,T;L^{2}(\Gamma ))}\leq C
$ and so}
\begin{equation*}
\left. \beta (\theta _{n})\right\vert_{{\Gamma} }\rightarrow \left. \beta
(\theta )\right\vert_{{\Gamma} }\mbox{ weakly in }L^{2}(0,T;L^{2}(\Gamma )),%
\mbox{ as }n\rightarrow \infty .
\end{equation*}
Passing to the limit as $n\rightarrow \infty $ in the weak form{s 
(cf.~\eqref{13}--\eqref{14})}
\begin{align*}
 \int_{0}^{T}\left\langle \frac{d\theta _{n}}{dt}(t),\psi
_{1}(t)\right\rangle _{V^{\prime },V}dt+\int_{Q}\nabla \beta (\theta
_{n})\cdot \nabla \psi _{1}\,{dx\hskip1pt dt}{}+\int_{Q}\frac{d\varphi
_{n}}{dt}\psi _{1}\,{dx\hskip1pt dt} \\
+\int_{\Sigma }\alpha \beta (\theta _{n})\psi _{1}\,{ds\hskip1pt dt}{}
=\int_{Q}u_{n}\psi _{1}\,{dx\hskip1pt dt}{}+\int_{\Sigma }\alpha
v_{n}\psi _{1}\,{ds\hskip1pt dt},
\end{align*}
\begin{align*}
 \int_{Q}\frac{d\varphi _{n}}{dt}\psi _{2}\,{dx\hskip1pt dt}
+\int_{Q}\nabla \varphi _{n}\cdot \nabla \psi _{2}\,{dx\hskip1pt dt}
+\int_{Q}(\varphi _{n}^{3}-\varphi _{n})\psi _{2}\,{dx\hskip1pt dt} \\
 =\int_{Q}\left( \frac{1}{\theta _{c}}-\frac{1}{\theta _{n}}\right) \psi
_{2}\,{dx\hskip1pt dt},
\end{align*}
{for all $\psi_1, \, \psi_2 \in L^2(0,T;V)$,} we obtain by the previous
convergences that $(\theta ,\varphi )$ satisfies (\ref{13})--(\ref{14}),
which means that it is the solution to (\ref{PF1})-(\ref{PF6}) corresponding
to $u$ and $v.$

Finally, we pass to the limit in (\ref{22}) using the weakly lower
semicontinuity property of the terms in $J$ and get
\begin{equation*}
J(u,v,\eta )=d.
\end{equation*}
This concludes the proof, by specifying that $(u,v,\eta )$ and the
corresponding states $(\theta ,\varphi )$ are optimal in $(P)$.\hfill $
\square $

\section{Approximating control problem}

In this section we {consider} an approximating problem $(P_{\varepsilon })$
and show its convergence in a {suitable} sense to $(P)$. First, we
{introduce the convex function
\begin{equation}
j:\mathbb{R}\rightarrow \mathbb{R},\quad j(r)=\left\vert r-\theta
_{c}\right\vert   \label{j}
\end{equation}
whose subdifferential is the graph $H$ defined in \eqref{defH}.
The conjugate of $j$ is}
\begin{equation*}
j^{\ast }(\omega )=\sup\limits_{r\in \mathbb{R}}\,(\omega r-j(r))
\end{equation*}
and {precisely} reads
\begin{equation}
j^{\ast }(\omega )=\omega \theta _{c}+I_{[-1,1]}(\omega ).  \label{j*}
\end{equation}
We mention that, if $K$ is a closed convex set, we denote by $
I_{K}$ its indicator function, which is defined by $I_{K}(r)=0$ if $r\in K$, $
I_{K}(r)=+\infty $ otherwise. 

{Let us} recall that two conjugate functions $j$ and $j^{\ast }$ satisfy the
relations (see{, e.g., \cite[p.~6]{vb-springer-2010})
\begin{equation}\label{subdiff}
j(r)+j^{\ast }(\omega ) \geq r\omega \ \mbox{ for all }r\mbox{, }\omega \in
\mathbb{R}; \quad
j(r)+j^{\ast }(\omega ) = r\omega \ \mbox{ iff }\omega \in \partial j(r),
\end{equation}
 where $\partial j$ denotes the subdifferential of $j$. In our special case,  \eqref{subdiff} reduces to
\begin{align}
j(r)+\omega \theta_c -\omega r\geq 0 \  \mbox{ for all } r\in\RR, \ \omega \in
[-1,1]; \qquad \notag \\
j(r)+\omega \theta_c -\omega r = 0 \ \mbox{ iff }\omega \in H(r). \label{fenchel}
\end{align}
Then, we let $\varepsilon >0$ and state the approximating problem as follows. 
Setting
\begin{eqnarray*}
J_{\varepsilon }(u,v,\eta ) =\lambda _{1}\int_{Q}(\theta -\theta
_{f})^2{dx\hskip1pt dt}+\lambda _{2}\int_{Q}(\varphi -\eta )^2 {dx\hskip1pt dt} \\
+\frac{1}{\varepsilon }\int_{Q}(j(\theta )+\eta \theta _{c}-\eta \theta ){dx\hskip1pt dt},
\end{eqnarray*}
we deal with the minimization problem
\begin{equation}
\mbox{Minimize }J_{\varepsilon }(u,v,\eta )\mbox{ for all }(u,v,\eta )\in
K_{1}\times K_{2}\times K_{[-1,1]},  \tag{$P_{\varepsilon }$}
\end{equation}
subject to (\ref{PF1})-(\ref{PF6}), where
\begin{equation}
K_{[-1,1]}=\{\eta \in L^{\infty }(Q)\, :\ \left\vert \eta (t,x)\right\vert \leq 1
\mbox{ a.e. }(t,x)\in Q\}.  \label{K-eta}
\end{equation}}
\medskip

\noindent \textbf{Proposition 3.1. }\textit{Let the assumptions} 
(\ref{10-0})\textit{\ and }(\ref{PF9})\textit{\ hold. Then }$
(P_{\varepsilon })$\textit{\ has at least one solution.}

\medskip

\noindent \textbf{Proof. }According to (\ref{fenchel}) and \eqref{K-eta}, 
{we have that 
$d_{\varepsilon
}=\inf\limits_{u,v,\eta }J_{\varepsilon }(u,v,\eta )\geq 0.$ Let} $
(u_{\varepsilon }^{n},v_{\varepsilon }^{n},\eta _{\varepsilon }^{n})_{n}$ be
a minimizing sequence for $J_{\varepsilon },$ {that is,}
\begin{equation}
d_{\varepsilon }\leq J_{\varepsilon }(u_{\varepsilon }^{n},v_{\varepsilon
}^{n},\eta _{\varepsilon }^{n})\leq d_{\varepsilon }+\frac{1}{n}\, .  \label{23}
\end{equation}
As in Theorem~2.3 we obtain that
{
\begin{align*}
u_{\varepsilon }^{n}\rightarrow u_{\varepsilon }^{\ast }\ \mbox{ {weakly*} in }
L^{\infty }(Q),\quad
v_{\varepsilon }^{n}\rightarrow v_{\varepsilon }^{\ast }\ \mbox{ {weakly*} in }L^{\infty }(\Sigma ), \qquad \\
\eta _{\varepsilon }^{n}\rightarrow \eta _{\varepsilon }^{\ast }\ \mbox{ {weakly*}
in }L^{\infty }(Q),\ \mbox{ as }n\rightarrow \infty ,
\end{align*}
} and $u_{\varepsilon }^{\ast }\in K_{1},$ $v_{\varepsilon
}^{\ast }\in K_{2},$ $\eta _{\varepsilon }^{\ast }\in {K_{[-1,1]}}.$
Also, for the corresponding state {we infer that}
\begin{align*}
\theta _{\varepsilon }^{n} \rightarrow \theta _{\varepsilon }^{\ast }\
\mbox{ weakly in }L^{2}(0,T;V)\cap W^{1,2}([0,T];V^{\prime })\qquad \\
\mbox{and strongly in }L^{2}(Q),\mbox{ as }n \rightarrow \infty ,
\end{align*}
\begin{align*}
\varphi _{\varepsilon }^{n}\rightarrow \varphi _{\varepsilon }^{\ast }\
\mbox{ weakly in }L^{2}(0,T;H^{2}(\Omega ))\cap W^{1,2}([0,T];L^{2}(\Omega )), 
\qquad \qquad\qquad\\
\mbox{{weakly*} in }L^{\infty }(0,T;V),\mbox{ and strongly in }
L^{2}(0,T;V),\mbox{ as }n\rightarrow \infty ,
\end{align*}
and
\begin{equation*}
\eta _{\varepsilon }^{n}\theta _{\varepsilon }^{n}\rightarrow \eta
_{\varepsilon }^{\ast }\theta _{\varepsilon }^{\ast }\ \mbox{ weakly in }
L^{2}(Q),\mbox{ as }n\rightarrow \infty .
\end{equation*}
In a similar way as proved in Theorem 2.3 we deduce {that}
\begin{equation*}
\frac{1}{\theta _{\varepsilon }^{n}}\rightarrow \frac{1}{\theta
_{\varepsilon }^{\ast }}\ \mbox{ weakly in }L^{2}(0,T;V),\mbox{ as }
n\rightarrow \infty ,
\end{equation*}
\begin{equation*}
(\varphi _{\varepsilon }^{n})^{3}\rightarrow (\varphi _{\varepsilon }^{\ast
})^{3}\ \mbox{ weakly in }L^{2}(Q),\mbox{ as }n\rightarrow \infty ,
\end{equation*}
\begin{equation*}
\beta (\theta _{\varepsilon }^{n})\rightarrow \beta (\theta _{\varepsilon
}^{\ast })\ \mbox{ weakly in }L^{2}(0,T;V),\mbox{ as }n\rightarrow \infty ,
\end{equation*}
\begin{equation*}
\left. \beta (\theta _{\varepsilon }^{n})\right\vert _{{\Gamma}}\rightarrow
\left. \beta (\theta _{\varepsilon }^{\ast })\right\vert _{{\Gamma}}\
\mbox{
weakly in }L^{2}(0,T;L^{2}(\Gamma )),\mbox{ as }n\rightarrow \infty
\end{equation*}
and {then show} that $(\theta _{\varepsilon }^{\ast },\varphi _{\varepsilon
}^{\ast })$ is a solution to (\ref{PF1})--(\ref{PF6}) corresponding to ($
u_{\varepsilon }^{\ast },v_{\varepsilon }^{\ast }).$

Next, since $j$ is {Lipschitz} continuous and $\theta _{\varepsilon
}^{n}$
{converges strongly to $\theta _{\varepsilon }^{\ast }$ in
$L^2(Q)$,} we have
\begin{equation*}
j(\theta _{\varepsilon }^{n})\rightarrow j(\theta _{\varepsilon }^{\ast })\ %
{\mbox{ strongly in } L^2(Q)},\mbox{ as }n\rightarrow \infty .
\end{equation*}
{When} passing to the limit in (\ref{23}),
in the third term of $J_{\varepsilon }$  {we exploit the weak lower
semicontinuity. Then,} we get $J_{\varepsilon }(u_{\varepsilon }^{\ast
},v_{\varepsilon }^{\ast },\eta _{\varepsilon }^{\ast })=d_{\varepsilon }$.
In conclusion,  $(u_{\varepsilon }^{\ast },v_{\varepsilon }^{\ast },\eta
_{\varepsilon }^{\ast })$ and the corresponding state $(\theta _{\varepsilon
}^{\ast },\varphi _{\varepsilon }^{\ast })$ are optimal in $(P_{\varepsilon
}).$\hfill $\square $

\medskip

The next theorem proves that $(P_{\varepsilon })$ converges to $(P)$ %
{in some sense} as $\varepsilon \rightarrow 0.$

\medskip

\noindent \textbf{Theorem 3.2. }\textit{Under the hypotheses of Theorem~%
{2.3}, {for any $\varepsilon >0$} let the pair} $\{(u_{\varepsilon
}^{\ast },v_{\varepsilon }^{\ast }, \eta_{\varepsilon }^{\ast
}),(\theta_{\varepsilon }^{\ast }, \varphi_{\varepsilon }^{\ast })\}$
\textit{be optimal in} $(P_{\varepsilon }).$\textit{\ Then,
{we have
that}}
\begin{align}
&u_{\varepsilon }^{\ast }\rightarrow u^{\ast }\
\mbox{
\textit{{weakly*} in} }L^{\infty }(Q),\mbox{ \textit{as} }
{\varepsilon \rightarrow 0} ,  \label{24} \\
&v_{\varepsilon }^{\ast }\rightarrow v^{\ast }\
\mbox{
\textit{{weakly*} in} }L^{\infty }(\Sigma ),\mbox{ \textit{as} }
{\varepsilon \rightarrow 0} ,  \label{25} \\
&\eta _{\varepsilon }^{\ast }\rightarrow \eta ^{\ast }
\mbox{ \textit{{weakly*} in}
}L^{\infty }(Q), \mbox{ \textit{as} }{\varepsilon \rightarrow 0} ,
\label{26}
\end{align}
\begin{align}
\theta _{\varepsilon }^{\ast }\rightarrow \theta ^{\ast }
\mbox{ \mbox{weakly
in} }L^{2}(0,T;V)\cap W^{1,2}([0,T];V^{\prime }) \qquad  \notag \\
\mbox{\textit{and strongly in} }L^{2}(Q),\mbox{
\textit{as} }{\varepsilon \rightarrow 0} ,  \label{27}
\end{align}
\begin{align}
\varphi _{\varepsilon }^{\ast }\rightarrow \varphi ^{\ast }
\mbox{
\textit{weakly in} }L^{2}(0,T;H^{2}(\Omega ))\cap W^{1,2}([0,T];L^{2}(\Omega
)), \qquad\qquad\qquad  \notag \\
\mbox{\textit{{weakly*} in} }L^{\infty }(0,T;V),%
\mbox{ \textit{and
strongly in} }L^{2}(0,T;V),\mbox{ \textit{as} }{\varepsilon \rightarrow
0} ,\label{28}
\end{align}
\textit{ where }$(\theta ^{\ast },\varphi ^{\ast })$\textit{\ is the
solution to} (\ref{PF1})--(\ref{PF6}) \textit{corresponding to }$(u^{\ast
},v^{\ast },\eta ^{\ast })$\textit{\ and the {pair} }$\{(u^{\ast
},v^{\ast },\eta ^{\ast }),(\theta ^{\ast },\varphi ^{\ast })\}$ \textit{is
optimal in }$(P).$
{\textit{Furthermore, every triplet $(u^{\ast},v^{\ast },\eta ^{\ast })$
obtained in \eqref{24}--\eqref{26} as weak* limits of subsequences of
$\{(u_{\varepsilon }^{\ast },v_{\varepsilon }^{\ast },
\eta_{\varepsilon }^{\ast })\}$ yields an optimal solution~to~$(P).$}}

\medskip

\noindent \textbf{Proof. }Let $\{(u_{\varepsilon }^{\ast },v_{\varepsilon
}^{\ast },\eta _{\varepsilon }^{\ast }),(\theta _{\varepsilon }^{\ast
},\varphi _{\varepsilon }^{\ast })\}$ be optimal in $(P_{\varepsilon }).$
Then we can write
\begin{equation*}
J_{\varepsilon }(u_{\varepsilon }^{\ast },v_{\varepsilon }^{\ast },\eta
_{\varepsilon }^{\ast })\leq J_{\varepsilon }(u,v,\eta ),
\end{equation*}
for any $(u,v,\eta )\in K_{1}\times K_{2}\times K_{[-1,1]}.$ In particular,
we set $u=\widetilde{u},$ $v=\widetilde{v},$ $\eta =\widetilde{\eta },$
where $(\widetilde{u},\widetilde{v},\widetilde{\eta })$ is a solution to ($
P) $ with the corresponding state $\widetilde{\theta },$ $\widetilde{\varphi
}$ {solving} (\ref{PF1})--(\ref{PF6}). This entails that $\widetilde{%
\eta }\in H(\, \widetilde{\theta } \, )
{{}\equiv \partial j
(\, \widetilde{\theta } \, )}$ {a.e. in} $Q.$ The previous inequality
reads
\begin{align}
\lambda _{1}\int_{Q}(\theta _{\varepsilon }^{\ast }-\theta _{f})^{2}\,%
{dx\hskip1pt dt} +\lambda _{2}\int_{Q}(\varphi _{\varepsilon }^{\ast
}-\eta _{\varepsilon }^{\ast })^{2}\,{dx\hskip1pt dt} +\frac{1}{%
\varepsilon }\int_{Q}(j(\theta _{\varepsilon }^{\ast })+\eta _{\varepsilon
}^{\ast }\theta _{c}-\eta _{\varepsilon }^{\ast }\theta _{\varepsilon
}^{\ast })\,{dx\hskip1pt dt}  \notag \\
\leq \lambda _{1}\int_{Q}(\widetilde{\theta }-\theta_{f})^{2}\,
{dx\hskip1pt dt} +\lambda _{2}\int_{Q}(\widetilde{\varphi }-\widetilde{%
\eta })^{2}\,{dx\hskip1pt dt} +\frac{1}{\varepsilon }\int_{Q}(j(%
\widetilde{\theta })+\widetilde{\eta }\theta _{c}-\widetilde{\eta }
\widetilde{\theta })\,{dx\hskip1pt dt} ,  \label{29}
\end{align}
and we see by (\ref{fenchel}) that the last term on the right-hand side is
actually zero. Then the right-hand side is bounded by a constant.

{In view of \eqref{K1}--\eqref{K2} and \eqref{K-eta}, by} the
boundedness of the optimal controllers $(u_{\varepsilon }^{\ast
},v_{\varepsilon }^{\ast },\eta _{\varepsilon }^{\ast })$ we {obtain} (%
\ref{24})-(\ref{26}).
{Thanks to Theorem~2.2 (cf. especially (\ref{15})--(\ref{16})), it is
straightforward to deduce} (\ref{27})--(\ref{28}). Then, writing the weak %
{formulations} (\ref{13})--(\ref{14}) for the approximating state and
passing to the limit as $\varepsilon \rightarrow 0$ we deduce that $(\theta
^{\ast },\varphi ^{\ast })$ is the solution to (\ref{PF1})--(\ref{PF6})
corresponding to $(u^{\ast },v^{\ast },\eta ^{\ast }).$

Finally, we have to show that $\eta ^{\ast }\in H(\theta ^{\ast })$ %
{a.e. in} $Q.$ We {set}
\begin{equation*}
\zeta _{\varepsilon }=\int_{Q}(j(\theta _{\varepsilon }^{\ast })
+ {}{ \eta _{\varepsilon }^{\ast } \theta_c}
-\eta _{\varepsilon }^{\ast }\theta
_{\varepsilon }^{\ast })\,{dx\hskip1pt dt}
\end{equation*}
{and remark that}
\begin{equation*}
0\leq {\gamma _{\varepsilon }:= \frac{1}{\varepsilon }\zeta
_{\varepsilon }\leq C}
\end{equation*}
{for some constant $C$  (independent of $\varepsilon$), because of \eqref{29}.  Hence,
we have that} $\zeta _{\varepsilon }=\varepsilon \gamma _{\varepsilon
}\rightarrow 0$, as $\varepsilon \rightarrow 0.$
{On the other hand,
passing to the limit we infer that}
\begin{equation*}
0\leq \int_{Q}(j(\theta ^{\ast })+\eta ^{\ast }\theta _{c}-\eta ^{\ast
}\theta ^{\ast })\,{dx\hskip1pt dt} \leq \lim_{\varepsilon \rightarrow
0}\zeta _{\varepsilon }=0.
\end{equation*}
We deduce that \, $j(\theta ^{\ast }(t,x))+ {}{ \eta^{\ast }(t,x)\theta_c}-\eta ^{\ast }(t,x)\theta ^{\ast }(t,x)=0$ \,
{a.e. $(t,x)\in
Q$}\, {and, thanks to \eqref{fenchel},} this implies 
that \, $\eta ^{\ast }(t,x)\in \partial j(\theta
^{\ast }(t,x))=H(\theta ^{\ast }(t,x))$ \, a.e. $(t,x)\in Q.$

{Then, we} pass to the limit in (\ref{29}) as $\varepsilon \rightarrow
0 $ and obtain%
\begin{equation*}
J(u^{\ast },v^{\ast },\eta ^{\ast })\leq J(\widetilde{u},\widetilde{v},%
\widetilde{\eta }),
\end{equation*}
for any $(\widetilde{u},\widetilde{v},\widetilde{\eta })\in K_{1}\times
K_{2}\times K_{3}$ {(cf.~\eqref{K3}),} with $(\widetilde{\theta },%
\widetilde{\varphi })$ solution to (\ref{PF1})--(\ref{PF6}). This shows that
$\{(u^{\ast },v^{\ast },\eta ^{\ast }),(\theta ^{\ast },\varphi ^{\ast })\}$
is optimal in\textit{\ }$(P).$\hfill $\square $

\section{Optimality conditions}

In this section we compute the optimality conditions for {the} problem $
(P_{\varepsilon }).$ We prove that whatever would be the optimal controllers
$u_{\varepsilon }^{\ast },$ $v_{\varepsilon }^{\ast },$ $\eta _{\varepsilon
}^{\ast },$ they are represented by the expressions given in Theorem 4.6. To
this end, we use some intermediate results proved for a second approximating
problem, introduced in order to regularize the function $j.$ We recall that
the Moreau{-Yosida} regularization is defined {by}
\begin{equation}
j_{\sigma }(r)=\inf_{s\in \mathbb{R}}\left\{ \frac{\left\vert r-s\right\vert
^{2}}{2\sigma }+j(s)\right\} ,\mbox{ for any }r\in \mathbb{R},\ \sigma >0,
\label{moreau}
\end{equation}
and that it can be still written {as}
\begin{equation}
j_{\sigma }(r)=\frac{1}{2\sigma }\left\vert (I+\sigma H)^{-1}r-r\right\vert
^{2}+j((I+\sigma H)^{-1}r),  \label{moreau-1}
\end{equation}
where $I$ is the identity on $\mathbb{R}.$ The function $j_{\sigma }$ is
convex, Lipschitz continuous {along with its derivative}, and {it}
has the properties (see {\cite[p.~48]{vb-springer-2010}}):%
\begin{equation}
{0\leq j_{\sigma }(r) \leq j(r) \ \mbox{ and } \  \lim_{\sigma
\rightarrow 0}j_{\sigma }(r) = j(r),\mbox{ for any }r\in \mathbb{R}. }
\label{30}
\end{equation}

Let $(u_{\varepsilon }^{\ast },v_{\varepsilon }^{\ast },\eta _{\varepsilon
}^{\ast })$ be optimal in $(P_{\varepsilon }).$ Following a
technique developed in \cite{vbpit}, we introduce the
approximating penalized problem:
\begin{equation}
\mbox{Minimize }J_{\varepsilon ,\sigma }(u,v,\eta )\mbox{ for all }(u,v,\eta )\in
K_{1}\times K_{2}\times K_{[-1,1]},  \tag{$P_{\varepsilon,\sigma }$}
\end{equation}subject to (\ref{PF1})--(\ref{PF6}), where
\begin{align*}
J_{\varepsilon ,\sigma }(u,v,\eta )=\lambda
_{1}\!\int_{Q}(\theta -\theta_{f})^{2}{dx\hskip1pt dt}
+\lambda _{2}\! \int_{Q}(\varphi -\eta )^{2}{dx\hskip1pt dt} 
+\frac{1}{\varepsilon }\int_{Q}(j_{\sigma }(\theta )+\eta \theta _{c}-\eta
\theta )\,{dx\hskip1pt dt} \\
 +\int_{Q}(u-u_{\varepsilon }^{\ast })^{2}\,{dx\hskip1pt dt}
 +\int_{\Sigma}(v-v_{\varepsilon }^{\ast })^{2}\,{ds\hskip1pt dt} 
 +\int_{Q}(\eta -\eta _{\varepsilon }^{\ast })^{2}\,{dx\hskip1pt dt},
\end{align*}
and $K_{[-1,1]}$ is defined by (\ref{K-eta}).

It is obvious that problem $(P_{\varepsilon ,\sigma })$ has at least one
solution, and the proof {can be done arguing as in} Proposition~3.1. %
{Now,} we shall show that {in a formal way} $(P_{\varepsilon
,\sigma })\rightarrow (P_{\varepsilon })$ as $\sigma \rightarrow 0.$

\medskip

\noindent \textbf{Proposition 4.1. }\textit{{Assume that 
\eqref{10-0} and \eqref{PF9} hold.} Let }$(u_{\varepsilon }^{\ast
},v_{\varepsilon }^{\ast },\eta _{\varepsilon }^{\ast })$\textit{\ be
optimal in }$(P_{\varepsilon })$\textit{\ and }$(u_{\varepsilon ,\sigma
}^{\ast },v_{\varepsilon ,\sigma }^{\ast },\eta _{\varepsilon ,\sigma
}^{\ast })$\textit{\ be optimal in }$(P_{\varepsilon ,\sigma }).$\textit{\
Then, {we have that}}
\begin{align}
& u_{\varepsilon ,\sigma }^{\ast }\rightarrow u_{\varepsilon }^{\ast }\
\mbox{
\textit{{weakly*} in} }L^{\infty }(Q)%
\mbox{ \textit{{and}
strongly in} }L^{2}(Q),\mbox{ \textit{as} }\sigma \rightarrow 0,
\label{pier2} \\
& v_{\varepsilon ,\sigma }^{\ast }\rightarrow v_{\varepsilon }^{\ast }\
\mbox{
\textit{{weakly*} in} }L^{\infty }(\Sigma )%
\mbox{ \textit{{and}
strongly in} }L^{2}(\Sigma ),\mbox{ \textit{as} }\sigma \rightarrow 0,
\label{31} \\
& \eta _{\varepsilon ,\sigma }^{\ast }\rightarrow \eta _{\varepsilon }^{\ast
}\ \mbox{ \textit{{weakly*} in} }L^{\infty }(Q)%
\mbox{ \textit{{and}  strongly in}
}L^{2}(Q),\mbox{ \textit{as} }\sigma \rightarrow 0,  \label{pier3}
\end{align}
\textit{and the corresponding states }$(\theta _{\varepsilon ,\sigma }^{\ast
},{\varphi}_{\varepsilon ,\sigma }^{\ast })$\textit{\ converge to }
\textit{the optimal states}\textit{\ }$(\theta
_{\varepsilon }^{\ast },\varphi _{\varepsilon }^{\ast })$\textit{\
{that
correspond to $(u_{\varepsilon }^{\ast
},v_{\varepsilon }^{\ast },\eta _{\varepsilon }^{\ast })$} in }$
(P_{\varepsilon }).$\textit{\ }

\medskip

\noindent \textbf{Proof. }We write that $(u_{\varepsilon ,\sigma }^{\ast},
v_{\varepsilon ,\sigma }^{\ast },\eta _{\varepsilon ,\sigma }^{\ast })$
\textit{\ }is optimal in\textit{\ }$(P_{\varepsilon ,\sigma }),$ that is%
\begin{align}
&\lambda _{1}\int_{Q}({\theta}^{\ast }_{\varepsilon ,\sigma}-\theta
_{f})^{2}\,{dx\hskip1pt dt} +\lambda _{2}\int_{Q}(\varphi _{\varepsilon
,\sigma }^{\ast }-\eta _{\varepsilon ,\sigma }^{\ast })^{2}\,%
{dx\hskip1pt dt}  \notag \\
&{}+ \frac{1}{\varepsilon }\int_{Q}(j_{\sigma }(\theta _{\varepsilon ,\sigma
}^{\ast })+\eta _{\varepsilon ,\sigma }^{\ast }\theta _{c}-\eta
_{\varepsilon ,\sigma }^{\ast }{\theta}^{\ast }_{\varepsilon ,\sigma})\,
{dx\hskip1pt dt}  \notag \\
&{}+\int_{Q}(u_{\varepsilon ,\sigma }^{\ast }-u_{\varepsilon }^{\ast })^{2}\,
{dx\hskip1pt dt} +\int_{\Sigma }(v_{\varepsilon ,\sigma }^{\ast
}-v_{\varepsilon }^{\ast })^{2}\,{ds\hskip1pt dt} +\int_{Q}(\eta
_{\varepsilon ,\sigma }^{\ast }-\eta _{\varepsilon }^{\ast })^{2}\,%
{dx\hskip1pt dt}  \notag \\
&{}\leq \lambda _{1}\int_{Q}(\theta -\theta_{f})^{2}\,{dx\hskip1pt dt}
+\lambda _{2}\int_{Q}(\varphi -\eta )^{2}\,{dx\hskip1pt dt} +\frac{1}{%
\varepsilon }\int_{Q}(j_{\sigma }(\theta )+\eta \theta _{c}-\eta \theta )\,
{dx\hskip1pt dt}  \notag \\
&{}\qquad +\int_{Q}(u-u_{\varepsilon }^{\ast })^{2}\,{dx\hskip1pt dt}
+\int_{\Sigma }(v-v_{\varepsilon }^{\ast })^{2}\,{ds\hskip1pt dt}
+\int_{Q}(\eta -\eta _{\varepsilon }^{\ast })^{2}\,{dx\hskip1pt dt} ,
\label{32}
\end{align}
for {all} $u\in K_{1},$ $v\in K_{2},$ $\eta \in K_{[-1,1]}$, with $
(\theta ,\varphi )$ {denoting} the corresponding solution to (\ref{PF1}
)--(\ref{PF6}).

In particular, we set $u=u_{\varepsilon }^{\ast },$ $
v=v_{\varepsilon }^{\ast },$ $\eta =\eta _{\varepsilon }^{\ast }$ in
\eqref{32}. This
leads us to consider the corresponding solutions $\theta =\theta
_{\varepsilon }^{\ast },$ $\varphi =\varphi _{\varepsilon }^{\ast }$ to
(\ref{PF1})--(\ref{PF6}) as well. It follows that the left-hand side in (%
\ref{32}) is bounded independently of $\sigma ,$ because on the right-hand
side the last three terms vanish and $j_{\sigma }(\theta _{\varepsilon
}^{\ast })\leq j(\theta _{\varepsilon }^{\ast })$
{a.e. in $Q$, thanks
to} (\ref{30}). Consequently, by selecting subsequences ({still denoted}
by the subscript {${\sigma }$}) we get
{
\begin{align}
u_{\varepsilon ,\sigma }^{\ast } \rightarrow  u_{\varepsilon }\ \mbox{ {weakly*} in }L^{\infty }(Q), \quad v_{\varepsilon ,\sigma }^{\ast } \rightarrow v_{\varepsilon }\
\mbox{ {weakly*} in }L^{\infty }(\Sigma ),\notag \qquad \\ \eta _{\varepsilon ,\sigma }^{\ast } \rightarrow \eta _{\varepsilon }\
\mbox{ {weakly*} in }L^{\infty }(Q), \mbox{ as }\sigma \rightarrow 0. \label{pier4}
\end{align}}Relying on the estimates (\ref{15})--(\ref{16}) for the state
system we have%
\begin{align}
{\theta}^{\ast }_{\varepsilon ,\sigma} \rightarrow \theta _{\varepsilon }
\mbox{ weakly in }L^{2}(0,T;V)\cap W^{1,2}([0,T];V^{\prime })\qquad   \notag
\\
\mbox{and strongly in }L^{2}(Q),\mbox{ as }\sigma  \rightarrow 0,
\label{pier8}
\end{align}
\begin{align}
{\varphi}_{\varepsilon ,\sigma }^{\ast } \rightarrow \varphi _{\varepsilon }
\mbox{ weakly in }L^{2}(0,T;H^{2}(\Omega ))\cap W^{1,2}([0,T];L^{2}(\Omega
)),\qquad \qquad   \notag \\
\mbox{{weakly*} in }L^{\infty }(0,T;V),\mbox{ and strongly in }
L^{2}(0,T;V),\mbox{ as }\sigma  \rightarrow 0,  \label{pier9}
\end{align}
{where} $(\theta _{\varepsilon },\varphi _{\varepsilon })$ is the
solution to (\ref{PF1})--(\ref{PF6}) corresponding to $(u_{\varepsilon
},v_{\varepsilon },{\eta_{\varepsilon }}).$

Next, we pass to the limit in (\ref{32}) as $\sigma \rightarrow 0.$ First,
we assert that
\begin{equation}
\int_{Q}j(\theta _{\varepsilon })\,{dx\hskip1pt dt} \leq
\liminf\limits_{\sigma \rightarrow 0}\int_{Q}j_{\sigma }(\theta
_{\varepsilon ,\sigma }^{\ast })\,{dx\hskip1pt dt} ,  \label{33}
\end{equation}
where $\theta _{\varepsilon }$
{is the limit of $\theta
_{\varepsilon ,\sigma }^{\ast }.$ Indeed, by} (\ref{moreau-1}) we have
\begin{equation*}
\frac{1}{2\sigma }\int_{Q}\left\vert (I+\sigma \partial j)^{-1}\theta
_{\varepsilon ,\sigma }^{\ast }-\theta _{\varepsilon ,\sigma }^{\ast
}\right\vert ^{2}\,{dx\hskip1pt dt} \leq \int_{Q}j_{\sigma }(\theta
_{\varepsilon ,\sigma }^{\ast })\,{dx\hskip1pt dt} \leq
\mbox{
constant,}
\end{equation*}
which implies that
\begin{equation*}
\lim_{\sigma \rightarrow 0}\left\Vert (I+\sigma \partial j)^{-1}\theta
_{\varepsilon ,\sigma }^{\ast }-\theta _{\varepsilon ,\sigma }^{\ast
}\right\Vert _{L^{2}(Q)}=0.
\end{equation*}
Therefore, we deduce that
\begin{equation}
{(I+\sigma \partial j)^{-1}{\theta}^{\ast }_{\varepsilon ,\sigma} \to
\theta _{\varepsilon }\ \mbox{ strongly in }L^{2}(Q) \, \hbox{ as } \sigma
\rightarrow 0 .}  \label{33-0}
\end{equation}
Next, again by (\ref{moreau-1}) we can infer that%
\begin{equation*}
{ \int_{Q}j(\theta _{\varepsilon })\,{dx\hskip1pt dt} =
{\lim \limits_{\sigma \rightarrow 0}}\int_{Q}j((I+\sigma \partial
j)^{-1}{\theta}^{\ast }_{\varepsilon ,\sigma})\,{dx\hskip1pt dt} \leq
\liminf\limits_{\sigma \rightarrow 0}\int_{Q}j_{\sigma }(\theta
_{\varepsilon ,\sigma }^{\ast })\,{dx\hskip1pt dt} }
\end{equation*}
by {the Lipschitz continuity of $j$} and (\ref{33-0}).
Then, passing to the limit in (\ref{32}) as $\sigma \rightarrow 0$ we get%
\begin{align*}
&\lambda _{1}\int_{Q}(\theta _{\varepsilon }-\theta_{f})^{2}\,
{dx\hskip1pt dt} +\lambda _{2}\int_{Q}(\varphi _{\varepsilon }-\eta
_{\varepsilon })^{2}\,{dx\hskip1pt dt} +\frac{1}{\varepsilon }
\int_{Q}(j(\theta _{\varepsilon })+\eta _{\varepsilon }\theta _{c}-\eta
_{\varepsilon }\theta _{\varepsilon })\,{dx\hskip1pt dt} \\
&{}+\int_{Q}(u_{\varepsilon }-u_{\varepsilon }^{\ast })^{2}\,
{dx\hskip1pt dt} +\int_{\Sigma }(v_{\varepsilon }-v_{\varepsilon
}^{\ast })^{2}\,{ds\hskip1pt dt} +\int_{Q}(\eta _{\varepsilon }-\eta
_{\varepsilon }^{\ast })^{2}\,{dx\hskip1pt dt} \\
&{}\leq \lambda _{1}\int_{Q}(\theta _{\varepsilon }^{\ast }-\theta
_{f})^{2}\,{dx\hskip1pt dt} +\lambda _{2}\int_{Q}(\varphi _{\varepsilon
}^{\ast }-\eta _{\varepsilon }^{\ast })^{2}\,{dx\hskip1pt dt} +\frac{1}{%
\varepsilon }\int_{Q}(j(\theta _{\varepsilon }^{\ast })+\eta _{\varepsilon
}^{\ast }\theta _{c}-\eta _{\varepsilon }^{\ast }\theta _{\varepsilon
}^{\ast })\,{dx\hskip1pt dt} \\
&{}\leq \lambda _{1}\int_{Q}(\theta _{\varepsilon }-\theta _{f})^{2}\,
{dx\hskip1pt dt} +\lambda _{2}\int_{Q}(\varphi _{\varepsilon }-\eta
_{\varepsilon })^{2}\,{dx\hskip1pt dt} +\frac{1}{\varepsilon }
\int_{Q}(j(\theta _{\varepsilon })+\eta _{\varepsilon }\theta _{c}-\eta
_{\varepsilon }\theta _{\varepsilon })\,{dx\hskip1pt dt} .
\end{align*}
The second inequality can be written because $(u_{\varepsilon }^{\ast
},v_{\varepsilon }^{\ast },\eta _{\varepsilon }^{\ast })$\textit{\ }is
optimal in $(P_{\varepsilon }).$ {Hence, it is not difficult to see}
that
\begin{equation*}
u_{\varepsilon }=u_{\varepsilon }^{\ast },\ v_{\varepsilon }=v_{\varepsilon
}^{\ast },\ \eta _{\varepsilon }=\eta _{\varepsilon }^{\ast }\
\mbox{ a.e.
in }Q
\end{equation*}
{and consequently} $\theta _{\varepsilon }=\theta _{\varepsilon }^{\ast
}$ and $\varphi _{\varepsilon }=\varphi _{\varepsilon }^{\ast }$
{a.e.
in} $Q.$ Actually, going back to (\ref{32}) it follows that
{the convergences in \eqref{pier4} hold for the whole sequences and moreover
\begin{align*}
u_{\varepsilon ,\sigma }^{\ast }\rightarrow u_{\varepsilon }^{\ast }\ \mbox{
strongly in }L^{2}(Q), \quad v_{\varepsilon ,\sigma }^{\ast }\rightarrow v_{\varepsilon }^{\ast }\ \mbox{ strongly in }L^{2}(\Sigma ), \qquad \\
\eta _{\varepsilon ,\sigma }^{\ast }\rightarrow \eta _{\varepsilon }^{\ast } \
\mbox{ strongly in }L^{2}(Q),\ \mbox{ as }\sigma \rightarrow 0.
\end{align*}}This ends the proof. \hfill $\square $

\subsection{Optimality conditions for $(P_{\protect\varepsilon ,\protect%
\sigma })$}

For a later use we begin by proving the well-posedness of the problem%
\begin{equation}
W_{t}-a(t,x)\Delta W+b(t,x)\Phi =\omega (t,x),\ \mbox{ in }Q,  \label{70}
\end{equation}
\begin{equation}
\Phi _{t}-\Delta \Phi +c(t,x)\Phi +d(t,x)W_{t}=g(t,x),\ \mbox{ in }Q,
\label{70-0}
\end{equation}
\begin{equation}
-\frac{\partial W}{\partial \nu }=\alpha (x)(W-\gamma (t,x)),\ \mbox{ on }
\Sigma ,  \label{70-1}
\end{equation}
\begin{equation}
\frac{\partial \Phi }{\partial \nu }=0,\ \mbox{ on }\Sigma ,  \label{70-2}
\end{equation}
\begin{equation}
W(0)=0,\ \Phi (0)=0,\mbox{ in }\Omega . \   \label{70-3}
\end{equation}

\medskip

\noindent \textbf{Proposition 4.2. }\textit{Let the following conditions}
{
\begin{equation}
a,b,c,d\in L^{\infty }(Q),\quad
0<a_{0}\leq a(t,x)\leq \left\Vert a\right\Vert_{L^{\infty }(Q)}=:\left\vert
a\right\vert _{\infty },\mbox{ \textit{a.e.} } (t,x)\in Q,
\label{60-0}
\end{equation}\begin{equation}
{\omega, g}  \in L^{2}(Q),\quad
\gamma \in W^{1,2}([0,T],L^{2}(\Gamma )),\quad \alpha
\hbox{ \textit{satisfies }\eqref{PF9}}
  \label{60-1}
\end{equation}}

\noindent\textit{hold. Then, {the} problem }(\ref{70})--(\ref{70-3})
\textit{has a unique solution {$(W,\Phi)$ with}}
\begin{eqnarray*}
&W \in {}{L^{\infty }(0,T;V)\cap W^{1,2}([0,T];L^{2}(\Omega ))}, \\%
[.1cm]
&\Phi \in {}{L^{2}(0,T;H^{2}(\Omega ))\cap L^{\infty }(0,T;V)\cap
W^{1,2}([0,T];L^{2}(\Omega ))}.
\end{eqnarray*}
\textit{If $\, \gamma \equiv 0\, $ {in addition, we have that}}
\begin{equation}
W\in L^{2}(0,T;H^{2}(\Omega )).  \label{60-2}
\end{equation}

\medskip

\noindent \textbf{Proof. }We use a fixed point argument. In (\ref{70})--(\ref%
{70-3}) let us fix $\overline{\Phi }\in L^{2}(Q)$ and consider the equations %
{and conditions}
\begin{equation}
\overline{W}_{t}-a(t,x)\Delta \overline{W}+b(t,x)\overline{\Phi }=\omega
(t,x),\ \mbox{ in }Q,  \label{71}
\end{equation}
\begin{equation}
-\frac{\partial \overline{W}}{\partial \nu }=\alpha (x)(\overline{W}-\gamma
(t,x)),\ \mbox{ on }\Sigma ,  \label{72}
\end{equation}
\begin{equation}
\overline{W}(0)=0,\ \mbox{ in }\Omega ,  \label{73}
\end{equation}

\begin{equation}
\Phi _{t}-\Delta \Phi +c(t,x)\overline{\Phi }+d(t,x)\overline{W}
_{t}=g(t,x),\ \mbox{ in }Q,  \label{74}
\end{equation}
\begin{equation}
\frac{\partial \Phi }{\partial \nu }=0,\ \mbox{ on }\Sigma ,  \label{75}
\end{equation}
\begin{equation}
\Phi (0)=0,\ \mbox{ in }\Omega .  \label{76}
\end{equation}

\smallskip \noindent We
{first solve \eqref{71}--\eqref{73} and find $\overline{W}$,
then we put $\overline{W}$ in \eqref{74} and solve \eqref{74}--\eqref{76}
by finding $\Phi$. Thus, we construct a mapping
\begin{equation}
\Psi :C([0,T];L^{2}(\Omega ))\rightarrow C([0,T];L^{2}(\Omega ))
\quad \hbox{such that} \quad  \Psi (\overline{\Phi })=\Phi .
\label{Psi}
\end{equation}
We are going to show that a suitable power of $\Psi$} is a contraction.

First of all, we {claim} that (\ref{71})--(\ref{73}) has a unique
solution
\begin{equation}
\overline{W}\in {}{L^{\infty }(0,T;V)\cap W^{1,2}([0,T];L^{2}(\Omega ))}
.  \label{76-1}
\end{equation}
{Let us outline} the argument, without writing in detail all
computations. By taking a partition of the interval $[0,T]$,
{setting
$t_{i}=ih,$ $i=1,\ldots,{N},$ with} $h=T/{N},$ we consider the
system of finite differences%
{
\begin{align}
&\frac{1}{a_{i}(x)}\frac{w_{i}-w_{i-1}}{h}-\Delta w_{i} = f_{i}(x),\ \mbox{ in }\Omega ,  \label{76-2} \\
& {}- \frac{\partial w_{i}}{\partial \nu } =\alpha (x)(w_{i}-\gamma _{i}(x) ),\ \mbox{ on }\Gamma ,   \label{pier5}
\end{align}for  $i=1,\ldots,{N},$ with $w_0=0$. Here, $a_i$ denotes the mean value of $a$ on the time interval $((i-1)h, (ih))$, and the same definition can be set for $f_i$
provided $f $ is interpreted as $(\omega - b \overline{\Phi })/a$. On the other hand, in view of \eqref{60-1}, $\gamma_i$ can be defined as $\gamma_i(x) = \gamma( t_i , x)$, a.e. $x\in \Gamma $.}

{Given $w_{i-1}\in L^2(\Omega)$, the system \eqref{76-2}--\eqref{pier5} has a unique variational solution $w_{i}\in H^1(\Omega)$ such that
$\Delta w_{i}$ lies in $L^2(\Omega)$ and the normal derivative $\frac{\partial w_{i}}{\partial \nu }$ is in $L^2(\Gamma)$. Then, thanks to well-known elliptic regularity results, the finite difference scheme \eqref{76-2}--\eqref{pier5}, $i=1,\ldots,{N}$,}
has a unique solution $(w_{1},\ldots,w_{{N}})\in X^{{N}}$ where
\begin{equation*}
X=\{z\in H^{3/2}(\Omega );\ \Delta z\in L^{2}(\Omega )\}.
\end{equation*}
{An} a priori estimate is obtained by testing (\ref{76-2}) by $
w_{i}-w_{i-1},$ and summing with respect to $i$, from 1 to $k\leq {N}.$ %
{Recalling that $w_0 =0$, we obtain}
\begin{align*}
{\frac{h}{\left\vert a\right\vert _{\infty }}}\sum_{i=1}^{k}\left\Vert
\frac{w_{i}-w_{i-1}}{h}\right\Vert _{L^{2}(\Omega )}^{2}+\frac{1}{2}
\int_{\Omega }\left\vert \nabla w_{k}\right\vert ^{2}dx+\frac{1}{2}
\int_{\Gamma }\alpha _{m}\left\vert w_{k}\right\vert ^{2}ds \\
\leq \sum_{i=1}^{k}h\int_{\Omega }f_{i}\frac{w_{i}-w_{i-1}}{h}\,
dx+\sum_{i=1}^{k}h\int_{\Gamma }\alpha \, \gamma _{i}\frac{w_{i}-w_{i-1}}{h}
\,ds.
\end{align*}
The last term on the right-hand side can be written {as}
\begin{equation*}
\sum_{i=1}^{k}h\int_{\Gamma }\alpha \, \gamma _{i}\frac{w_{i}-w_{i-1}}{h}\,
ds=\int_{\Gamma }\alpha \, \gamma _{k}w_{k}\, ds-\sum_{i=1}^{k}\int_{\Gamma
}\alpha (\gamma _{i}-\gamma _{i-1})w_{i-1}\,ds
\end{equation*}
and standard computations involving the Young inequality and the discrete
Gronwall's lemma along with assumptions \eqref{60-0}--\eqref{60-1} lead to %
{some} estimates for the {functions} $\widehat{w}_{h}$ (%
{piecewise} linear in time interpolant) and {$\overline{w}_{h}$}
(piecewise constant {in time} interpolant). {Namely, we have}
\begin{equation*}
\left\Vert \partial _{t}\widehat{w}_{h}\right\Vert _{L^{2}(0,T;L^{2}(\Omega
))}+\left\Vert {\overline{w}_{h}}\right\Vert _{L^{\infty
}(0,T;H^{1}(\Omega ))}\leq C.
\end{equation*}
Using this we can pass to the limit (by {weak and weak*} compactness)
in the equations
\begin{eqnarray*}
&\displaystyle \frac{1}{\overline{a}_{h}}\partial _{t}\widehat{w}_{h}-\Delta
{}{\overline{w}_{h} }{}= \overline{f}_{h}, \ \mbox{ \ in }Q, & \\
\noalign{\smallskip} &-\displaystyle \frac{\partial \overline{w}_{h}}{%
\partial \nu } = {}{\alpha}{} (\overline{w}_{h}-\overline{\gamma }
_{h}), \ \mbox{ on }\Sigma , & \\
\noalign{\smallskip} &\widehat{w}_{h}(0,\, \cdot \, ) =0, \ \mbox{ in }
\Omega .&
\end{eqnarray*}
By comparison we also find the additional regularity
\begin{equation*}
\left\Vert \Delta \overline{w}_{h}\right\Vert _{L^{2}(0,T;L^{2}(\Omega
))}+\left\Vert \frac{\partial \overline{w}_{h}}{\partial \nu } \right\Vert
_{L^{2}(0,T;L^{2}(\Gamma ))}\leq C.
\end{equation*}
{We point out that if $\gamma \equiv 0,$ in the system \eqref{76-2}--\eqref{pier5} one can recover that the normal derivative of the solution on the boundary, i.e. $\frac{\partial w_{i}}{\partial \nu }$, belongs to $H^{1/2} (\Gamma) $, and consequently $w_i \in H^2(\Omega)$, whenever the product $\alpha \, w_{i}$ lies in $H^{1/2} (\Gamma) $. Now, we have that $\alpha$ satisfies \eqref{PF9}  and  the trace $w_{i}$ is in $H^{1} (\Gamma) $, on account of $w_i \in X$. Well, it is easy to check that the product of two elements of $H^{1} (\Gamma) $ belongs to
$W^{1,p} (\Gamma)$ for all $1\leq p<\infty$ due to the Leibniz rule and to the fact that $\Gamma$ is the (two-dimensional) boundary of a three-dimensional
domain $\Omega$. Hence, thanks to the $2D$ Sobolev embedding
$W^{1,p} (\Gamma) \subset  H^{1/2} (\Gamma) $ if $p \geq 4/3$, it turns out that
$$\frac{\partial w_{i}}{\partial \nu }= \alpha \, w_i \in H^{1/2} (\Gamma) .$$
Then, if $\gamma \equiv 0$, it is not difficult to obtain (\ref{60-2}) for $\overline{W}$.}

Consequently to (\ref{76-1}) it follows that
{the subsequent
linear parabolic problem} (\ref{74})--(\ref{76}) has a unique solution
\begin{equation*}
\Phi \in {}{L^{2}(0,T;H^{2}(\Omega ))\cap L^{\infty }(0,T;V)\cap
W^{1,2}([0,T];L^{2}(\Omega ))}. 
\end{equation*}

{At this point, we can write (\ref{71})--(\ref{73}) for two functions
$\overline{\Phi }_{1},$ $\overline{\Phi }_{2}\in L^{2}(Q)$ getting the respective solutions $\overline{W}_{1},$ $\overline{W}_{2}$ which satisfy (\ref{76-1}).
Then, take the difference,} divide by $a,$ multiply {the result} by ($
\overline{W}_{1}-\overline{W}_{2})_{t}$ and integrate over $Q.$
{After
a few standard computations, we obtain}
\begin{align}
\int_{0}^{t}\left\Vert {(\overline{W}_{1}-\overline{W}_{2})_{t}(\tau) }
\right \Vert _{L^{2}(\Omega )}^{2}d\tau +\left\Vert \nabla (\overline{W}
_{1}- \overline{W}_{2})(t)\right\Vert _{L^{2}(\Omega )}^{2}  \notag \\[1pt]
+\int_{\Gamma }\left|(\overline{W}_{1}-\overline{W}_{2})(t)\right|^{2}ds
\leq C_{1}\int_{0}^{t}\left\Vert (\overline{\Phi }_{1}-\overline{\Phi }
_{2})(\tau)\right\Vert _{L^{2}(\Omega )}^{2}d\tau  \label{77}
\end{align}
{for some constant $C_1$ depending only on the data in \eqref{60-0}--\eqref{60-1}. Next, we write (\ref{74})--(\ref{76})
for $\overline{\Phi }_{i},$ $\overline{W}_{i},$ $i=1,2$, subtract,
and test by ($\Phi _{1}-\Phi _{2})$, the difference of solutions. Hence,
it is a standard matter to infer that\begin{align}
&\left\Vert (\Phi _{1}-\Phi _{2})(t)\right\Vert _{L^{2}(\Omega )}^{2} \notag \\
&\leq
C_{2}\left(\int_{0}^{t}\left\Vert (\overline{\Phi }_{1}-\overline{\Phi }_{2})(\tau)\right\Vert _{L^{2}(\Omega )}^{2}d\tau +
\int_{0}^{t}\left\Vert {(\overline{W}_{1}-\overline{W}_{2})_{t}(\tau) }\right
\Vert _{L^{2}(\Omega )}^{2}d\tau \right) \label{pier6}
\end{align}
for some positive constant $C_2$. Then, combining \eqref{pier6} with \eqref{77}, we obtain
\begin{align}
\left\Vert (\Phi _{1}-\Phi _{2})(t)\right\Vert _{L^{2}(\Omega )}^{2}\leq
C_{3}\int_{0}^{t}\left\Vert (\overline{\Phi }_{1}-\overline{\Phi }_{2})(\tau)\right\Vert _{L^{2}(\Omega )}^{2}d\tau \quad \hbox{for all } t\in [0,T].  \label{79}
\end{align}
By observing that (cf.~\eqref{Psi}) $\Phi_i = \Psi(\overline{\Phi }_{i})$, $i=1,2$, is not difficult to check that relation (\ref{79}) implies
by recurrence that
\begin{equation}
\left\Vert \Psi^k(\overline{\Phi }_{1}) - \Psi^k(\overline{\Phi }_{2})
\right\Vert _{C([0,T];L^{2}(\Omega ))}^{2} \leq C_3 \frac{T^k}{k!}
\left\Vert \overline{\Phi }_{1} - \overline{\Phi }_{2}
\right\Vert _{C([0,T];L^{2}(\Omega ))}^{2}
\notag
\end{equation}
for all $k\in \NN.$ Then for $k$ large enough the above coefficient $C_3 \frac{T^k}{k!}$ is less than 1, and so $\Psi^k $ has a
unique fixed point  $\Phi $ which also fulfils $\Phi=\Psi (\Phi).$} \hfill $
\square $

\medskip

{In the sequel, we will assume the further regularity conditions \eqref{17} and \eqref{20} (besides \eqref{10}) in order we can take advantage of uniform $L^\infty$ estimates for both components of a solution to \eqref{PF1}--\eqref{PF6}.
Let us} resume the computation of the optimality conditions for $
(P_{\varepsilon ,\sigma }).$ Let $(u_{\varepsilon ,\sigma }^{\ast
},v_{\varepsilon ,\sigma }^{\ast },\eta _{\varepsilon ,\sigma }^{\ast })$
\textit{\ }and\textit{\ (}${\theta}^{\ast }_{\varepsilon ,\sigma},\varphi
_{\varepsilon ,\sigma }^{\ast })$ be optimal in\textit{\ }$(P_{\varepsilon
,\sigma })$ and $\lambda \in (0,1).$ We introduce the variations
\begin{equation*}
u_{\varepsilon ,\sigma }^{\lambda }=(1-\lambda )u_{\varepsilon ,\sigma
}^{\ast }+\lambda u=u_{\varepsilon ,\sigma }^{\ast }+\lambda \widetilde{u},\
\ { u \mbox{ arbitrary in }}K_{1},
\end{equation*}
\begin{equation*}
v_{\varepsilon ,\sigma }^{\lambda }=(1-\lambda )v_{\varepsilon ,\sigma
}^{\ast }+\lambda v=v_{\varepsilon ,\sigma }^{\ast }+\lambda \widetilde{v},\
\ { v \mbox{ arbitrary in }}K_{2},
\end{equation*}
\begin{equation*}
{\eta}_{\varepsilon ,\sigma }^{\lambda }=(1-\lambda )\eta _{\varepsilon
,\sigma }^{\ast }+\lambda \eta =\eta _{\varepsilon ,\sigma }^{\ast }+\lambda
\widetilde{\eta },\ \ { \eta \mbox{ arbitrary in }}K_{[-1,1]},
\end{equation*}
{with}
\begin{equation}
\widetilde{u}=u-u_{\varepsilon ,\sigma }^{\ast },\ \ \widetilde{v}
=v-v_{\varepsilon ,\sigma }^{\ast },\ \ \widetilde{\eta }=\eta -\eta
_{\varepsilon ,\sigma }^{\ast }\,.  \label{34}
\end{equation}
First of all, we note that the system (\ref{PF1})-(\ref{PF6})
corresponding to $u_{\varepsilon ,\sigma }^{\lambda }$ and $v_{\varepsilon
,\sigma }^{\lambda }$ has a unique solution $
{(\theta
_{\varepsilon ,\sigma }^{\lambda },\varphi _{\varepsilon ,\sigma
}^{\lambda })}$, and { 
\begin{equation}
\theta_{\varepsilon ,\sigma }^{\lambda }\rightarrow
{\theta}^{\ast }_{\varepsilon ,\sigma}, \quad \varphi _{\varepsilon ,\sigma
}^{\lambda }\rightarrow {\varphi}_{\varepsilon ,\sigma }^{\ast } \ \hbox{ strongly in }
L^{2}(Q), \, \hbox{ as } \lambda \rightarrow 0. \label{pier15}
\end{equation}
}
This can be obtained {by the estimate \eqref{cd} combined with weak* compactness. We~set 
\begin{equation}
\widetilde{\Theta }^{\lambda }=\frac{\theta _{\varepsilon ,\sigma }^{\lambda
}-{\theta}^{\ast }_{\varepsilon ,\sigma}}{\lambda },\mbox{ \ \ }\widetilde{
\Phi }^{\lambda }=\frac{\varphi _{\varepsilon ,\sigma }^{\lambda }-\varphi
_{\varepsilon ,\sigma }^{\ast }}{\lambda }  \label{34-1}
\end{equation}
and claim that 
\begin{equation}
\widetilde{\Theta }^{\lambda } \to Y \ \hbox{ and } \ 
\widetilde{\Phi }^{\lambda } \to \Phi  \ \mbox{ weakly in } L^2 (Q), \hbox{ as } 
\lambda \rightarrow 0, \label{pier21}
\end{equation}
where the limits $Y$ and $\Phi $ solve the system in variations}
\begin{equation}
Y_{t}-\Delta (\beta ^{\prime }(\theta _{\varepsilon ,\sigma }^{\ast
})Y)+\Phi _{t}=\widetilde{u},\ \mbox{ in }Q,  \label{35}
\end{equation}
\begin{equation}
\Phi _{t}-\Delta \Phi +(3(\varphi _{\varepsilon ,\sigma }^{\ast
})^{2}-1)\Phi =\frac{1}{({\theta}^{\ast }_{\varepsilon ,\sigma})^{2}}Y,\
\mbox{
in }Q,  \label{36}
\end{equation}
\begin{equation}
{}-\frac{\partial }{\partial \nu }(\beta ^{\prime }(\theta _{\varepsilon
,\sigma }^{\ast })Y)=\alpha (x)(\beta ^{\prime }(\theta _{\varepsilon
,\sigma }^{\ast })Y-\widetilde{v}),\ \mbox{ on }\Sigma ,  \label{37}
\end{equation}
\begin{equation}
\frac{\partial \Phi }{\partial \nu }=0,\ \mbox{ on }\Sigma ,  \label{38}
\end{equation}
\begin{equation}
{ Y(0)=0,\ \Phi (0)=0,\ \mbox{ in }\Omega . }  \label{40}
\end{equation}

The proof of {\eqref{pier21}} is done in Proposition~4.3, below.
Before that, we define a (very weak) \textit{solution} to (\ref{35})--(\ref{40}) as
a pair of functions $Y\in L^{2}(Q),$ $\Phi \in L^{2}(0,T;V)$ which satisfies
the system
\begin{align}
-\int_{Q}Y\psi_{t}{\hskip1pt dx\hskip1pt dt}-\int_{Q}\beta ^{\prime }(\theta _{\varepsilon
,\sigma }^{\ast })Y\Delta \psi {\hskip1pt dx\hskip1pt dt}-\int_{Q}\Phi \psi _{t}{\hskip1pt dx\hskip1pt dt}  \notag \\ \label{40-6}
=\int_{Q}\widetilde{u}\psi {\hskip1pt dx\hskip1pt dt}+\int_{\Sigma }\alpha \widetilde{v}\psi {\hskip1pt ds\hskip1pt dt}, 
\end{align}
\begin{align}
-\int_{Q}\Phi (\psi _{1})_{t}{\hskip1pt dx\hskip1pt dt}+\int_{Q}\nabla \Phi \cdot \nabla \psi_{1}{\hskip1pt dx\hskip1pt dt}+\int_{Q}(3({\varphi}_{\varepsilon ,\sigma }^{\ast })^{2}-1)\Phi \psi_{1}{\hskip1pt dx\hskip1pt dt} \notag  \\
=\int_{Q}\frac{1}{({\theta}^{\ast }_{\varepsilon ,\sigma})^{2}}Y\psi
_{1}{\hskip1pt dx\hskip1pt dt},  \label{pier26}
\end{align}
for all $\psi \in L^{2}(0,T;H^2(\Omega)) \cap W^{1,2}([0,T];L^{2}(\Omega ))$ 
solving the problem
\begin{align}
\psi_{t} + \Delta \psi = { - f_{Q}},  \ \mbox{ in } Q;  \qquad
\frac{\partial \psi }{\partial \nu } + \alpha \psi = 0, \ \mbox{ on }\Sigma ;
\label{40-7} \qquad
\psi (T) = 0 , \ \mbox{ in } \Omega, 
\end{align}
{for a generic $f_{Q}\in L^{2}(Q),$ and for all $\psi_1 \in L^{2}(0,T;V)\cap W^{1,2}([0,T];L^{2}(\Omega )) $ such that} $\psi_{1}(T)=0.$

\medskip

\noindent \textbf{Proposition 4.3. }\textit{Assume }(\ref{10-0}), (\ref{PF9}
), (\ref{17})\textit{\ and }(\ref{20})\textit{. }
\textit{Then the problem} (\ref{35})--(\ref{40}) \textit{has a unique
solution {$(Y, \Phi)$ with}}
\begin{equation}
Y\in L^{2}(0,T;L^{2}(\Omega )),\   \label{40-0}
\end{equation}
\begin{equation}
\Phi \in L^{2}(0,T;H^{2}(\Omega ))\cap L^{\infty }(0,T;V)\cap
W^{1,2}([0,T];L^{2}(\Omega )),  \label{40-1}
\end{equation}
\textit{and {the convergence properties in \eqref{pier21} hold.}}

\medskip
\noindent \textbf{Proof. }First, we prove the existence and
uniqueness of the solution to (\ref{35})--(\ref{40}). Due to the
hypotheses (\ref{20}){ ,  we infer} that the state system (\ref{PF1})--(%
\ref{PF6})%
{,
written for $u= u_{\varepsilon ,\sigma }^{\ast }$ and $v = v_{\varepsilon
,\sigma }^{\ast }$,} has the solution
{$ (\theta _{\varepsilon,\sigma }^{\ast }, {\varphi}_{\varepsilon ,\sigma }^{\ast })$ with both $\theta _{\varepsilon,\sigma }^{\ast }$ and $1 / \theta _{\varepsilon
,\sigma }^{\ast }$ bounded in $L^{\infty }(Q)$ (see (\ref{21-1}))
and consequently (cf.~\eqref{PF8}) $
\beta ^{\prime }({\theta}^{\ast }_{\varepsilon ,\sigma})\in L^{\infty }(Q).$
Moreover, in view of \eqref{17}, by (\ref{19}) we deduce the boundedness 
of ${\varphi}_{\varepsilon ,\sigma }^{\ast }$ in
$L^{\infty }(Q)\cap L^{\infty }(0,T;H^{2}(\Omega )).$}

We integrate (\ref{35}) and (\ref{37}) with respect to $\tau $ on $(0,t)$. %
{We obtain}
\begin{align*}
& Y(t,x)-\Delta \int_{0}^{t}\beta ^{\prime }(\theta _{\varepsilon ,\sigma
}^{\ast }(\tau ,x))Y(\tau ,x)d\tau +\Phi (t,x)=\int_{0}^{t}\widetilde{u}
(\tau ,x)d\tau , \\
& -\frac{\partial }{\partial \nu }\int_{0}^{t}(\beta ^{\prime }(\theta
_{\varepsilon ,\sigma }^{\ast }(\tau ,x))Y(\tau ,x))d\tau =\alpha
(x)\int_{0}^{t}(\beta ^{\prime }({\theta}^{\ast }_{\varepsilon ,\sigma}(\tau
,x))Y{(\tau ,x)}-\widetilde{v}(\tau ,x))d\tau 
\end{align*}
{and then set
$
W(t,x)=\int_{0}^{t}\beta ^{\prime }(\theta _{\varepsilon ,\sigma }^{\ast
}(\tau ,x))Y(\tau ,x)d\tau   $
for $(t,x) \in Q$, so that}
\begin{equation}
W_{t}(t,x)=\beta ^{\prime }(\theta _{\varepsilon ,\sigma }^{\ast
}(t,x))Y(t,x), \quad {(t,x)\in Q}.  \label{40-3}
\end{equation}
{Now, the system (\ref{35})--(\ref{40}) can be replaced} by
\begin{equation}
\Big(\,\frac{1}{\beta ^{\prime }({\theta}^{\ast }_{\varepsilon ,\sigma})}
W_{t}-\Delta W+\Phi \Big)(t,x)=\int_{0}^{t}\widetilde{u}(\tau ,x)d\tau ,\ \
(t,x)\in Q,  \label{40-4}
\end{equation}
\begin{equation}
\Phi_{t}-\Delta \Phi +(3{\varphi}_{\varepsilon ,\sigma }^{\ast }-1)\Phi
=\frac{1}{({\theta}^{\ast }_{\varepsilon ,\sigma})^{2}\beta ^{\prime
}({\theta}^{\ast }_{\varepsilon ,\sigma})}W_{t},\ \mbox{ in }Q, \label{pier23}
\end{equation}
\begin{equation}
-\frac{\partial W}{\partial \nu }(t,x)=\alpha (x)\Big(W(t,x)-\int_{0}^{t}
\widetilde{v}(\tau ,x) d\tau \Big),\ \ (t,x)\in \Sigma , \label{pier24}
\end{equation}
\begin{equation}
\frac{\partial \Phi }{\partial \nu }=0,\ \mbox{ on }\Sigma , \label{pier25}
\end{equation}
\begin{equation}
{ W(0)=0,\ \Phi (0)=0,\ \mbox{ in }\Omega . }  \label{40-5}
\end{equation}
{Here, we are allowed to apply Proposition~4.2 for
\begin{align*}
&a(t,x) = b(t,x)=\beta ^{\prime }(\theta _{\varepsilon }^{\ast }(t,x)),\mbox{
}\quad \omega (t,x)=\beta ^{\prime }(\theta _{\varepsilon }^{\ast
}(t,x))\int_{0}^{t}\widetilde{u}(\tau ,x)d\tau , \\
&c(t,x) =(3{\varphi}_{\varepsilon ,\sigma }^{\ast }(t,x)-1),\quad d(t,x)=-\frac{1}{({\theta}^{\ast }_{\varepsilon ,\sigma})^{2}\beta ^{\prime }(\theta
_{\varepsilon ,\sigma }^{\ast })},\quad g(t,x)=0,
\end{align*}
$(t, x)\in Q$, and
$\gamma (t,x) =\displaystyle \int_{0}^{t}\widetilde{v}(\tau ,x)d\tau$, $(t, x) \in \Sigma.$
Observing that such coefficients satisfy (\ref{60-0})--(\ref{60-1}),
we conclude that (\ref{40-4})--(\ref{40-5}) has a unique solution $(W,\Phi)$
with $W$ in the space $ L^{\infty }(0,T;V) \cap W^{1,2}([0,T];L^{2}(\Omega ))$.
Then, owing to (\ref{40-3}), it turns out that (\ref{40-0}) holds.}

{Next, we have to show that if $(Y, \Phi)$ fulfils the variational equalities in \eqref{40-6} and \eqref{pier26}, then  the pair $(W, \Phi)$
 with $W$ specified by \eqref{40-3} just solves the system (\ref{40-4})--(\ref{40-5}). Indeed, taking an arbitrary 
\begin{equation}
\zeta \in H^2(\Omega) \  \hbox{ such that } \ \frac{\partial \zeta }{\partial \nu } + \alpha \zeta = 0, \ \mbox{ on }\Gamma , \label{pier22}
\end{equation}
according to \eqref{40-7} we can choose $\psi (t,x) = (T-t) \zeta(x) $, $(x,t)\in Q$, in \eqref{40-6}. Then, if  $\zeta $ also belongs to ${\cal D} (\Omega)$, integrating by parts in time it is not difficult to recover the equality \eqref{40-4} in the sense of distributions in $\Omega$, for a.e. $t\in (0,T)$. Once \eqref{40-4} is proved, we can compare the terms and find additional regularity for $W$ (in particular, $\Delta W \in 
L^2(0,T;L^{2}(\Omega )$) in order to be able to get back to \eqref{40-6} and this time still use $\psi (t,x) = (T-t) \zeta(x) $, but with an auxiliary function $\zeta$ as in \eqref{pier22} to find the boundary condition \eqref{pier24} as well. 
A similar approach can be used on \eqref{pier26} taking now $\psi_1 (t,x) = (T-t) \zeta_1(x) $, $(x,t)\in Q$, with $\zeta_1$ arbitrary first in $H^1_0 (\Omega)$, then in 
$H^1 (\Omega)$ in order to arrive at an integrated version of \eqref{pier23}
and \eqref{pier25}. Then, it suffices to examine the regularity of $\Phi$ and realize that \eqref{40-1}, as well as \eqref{pier23}
and \eqref{pier25} directly, are satisfied.}

We prove now (\ref{pier21}). As mentioned in {Theorem~}2.2, the solution to (\ref{PF1})--(\ref{PF6}) is Lipschitz continuous with respect to the data.
Relying on (\ref{cd}) {and recalling \eqref{34-1},} we can write
\begin{eqnarray}
{\left\Vert \widetilde{\Theta }^{\lambda }\right\Vert
_{L^{2}(Q)}^{2}+\left\Vert \widetilde{\Phi }^{\lambda }\right\Vert
_{C([0,T];L^{2}(\Omega ))}^{2}+\left\Vert \widetilde{\Phi }^{\lambda }\right\Vert
_{L^{2}(0,T;V)}^{2}  \label{40-9} 
\leq \left\Vert \widetilde{u}\right\Vert _{L^{2}(Q)}^{2}+\left\Vert
\widetilde{v}\right\Vert _{L^{2}(\Sigma )}^{2}. } 
\end{eqnarray}
It is also obvious that, for each $\lambda ,$ the functions $\widetilde{%
\Theta }^{\lambda }$ and $\widetilde{\Phi }^{\lambda }$ are in the same
spaces as ${\theta}^{\ast }_{\varepsilon ,\sigma}$ and $\varphi
_{\varepsilon ,\sigma }^{\ast }$ are$,$ given by {Theorem~2.2}, and that they
satisfy the system
\begin{equation}
\widetilde{\Theta }_{t}^{\lambda }-\Delta \frac{\beta (\theta _{\varepsilon
,\sigma }^{\lambda })-\beta ({\theta}^{\ast }_{\varepsilon ,\sigma})}{%
\lambda }+\widetilde{\Phi }_{t}^{\lambda }=\widetilde{u},\ \mbox{ in }Q,
\label{40-10}
\end{equation}
\begin{equation}
\widetilde{\Phi }_{t}^{\lambda }-\Delta \widetilde{\Phi }^{\lambda
}+\left( (\varphi _{\varepsilon ,\sigma }^{\lambda })^{2}+\varphi _{\varepsilon
,\sigma }^{\lambda }{\varphi}_{\varepsilon ,\sigma }^{\ast }+(\varphi
_{\varepsilon ,\sigma }^{\ast })^{2}-1\right)\widetilde{\Phi }^{\lambda }=\frac{\widetilde{\Theta }^{\lambda }}{\theta _{\varepsilon ,\sigma }^{\lambda
}{\theta}^{\ast }_{\varepsilon ,\sigma}}
   ,\ \mbox{
in }Q, \label{pier27}
\end{equation}
\begin{equation}
{}-\frac{\partial }{\partial \nu }\frac{\beta (\theta _{\varepsilon ,\sigma
}^{\lambda })-\beta ({\theta}^{\ast }_{\varepsilon ,\sigma})}{\lambda }
=\alpha (x)\left( \frac{\beta (\theta _{\varepsilon ,\sigma }^{\lambda
})-\beta ({\theta}^{\ast }_{\varepsilon ,\sigma})}{\lambda }-\widetilde{v}
\right) ,\ \mbox{ on }\Sigma , \label{pier28}
\end{equation}
\begin{equation}
\frac{\partial \widetilde{\Phi }^{\lambda }}{\partial \nu }=0,\ \mbox{ on }
\Sigma , \label{pier29}
\end{equation}
\begin{equation}
\widetilde{\Theta}^{\lambda }  =0, \quad \widetilde{\Phi }^{\lambda }=0,
\ \mbox{ in } \Omega .  \label{40-11}
\end{equation}
{Thanks to (\ref{40-9}), at least for a subsequence of $\lambda$ we have} that
\begin{align*}
&\widetilde{\Theta }^{\lambda } \rightarrow \widetilde{Y} \ \mbox{
weakly in }L^{2}(Q),\, \mbox{ as }\lambda \rightarrow 0, \\
&\widetilde{\Phi }^{\lambda } \rightarrow \widetilde{\Phi } \ \mbox{ weakly in
}L^{2}(0,T;V)\mbox{ and weakly* in }L^{\infty }(0,T;L^{2}(\Omega )), \, \mbox{
as }\lambda \rightarrow 0.
\end{align*}
Let us test {(\ref{40-10}) by $\psi $ given as in (\ref{40-7}) and \eqref{pier27} 
by $\psi _{1}\in W^{1,2}([0,T];L^{2}(\Omega ))\cap L^{2}(0,T;V)$ with $\psi
_{1}(T)=0.$ Using the boundary conditions \eqref{pier28} and \eqref{pier29}, we} 
obtain%
\begin{align}
&-\int_{Q}\widetilde{\Theta }^{\lambda }\psi _{t}{\hskip1pt dx\hskip1pt dt}-\int_{Q}\frac{\beta
(\theta _{\varepsilon ,\sigma }^{\lambda })-\beta (\theta _{\varepsilon
,\sigma }^{\ast })}{\lambda }\Delta \psi {\hskip1pt dx\hskip1pt dt}
\notag \\
&-\int_{Q}\widetilde{\Phi}^{\lambda }\psi _{t}{\hskip1pt dx\hskip1pt dt}  
=\int_{Q}\widetilde{u}\psi {\hskip1pt dx\hskip1pt dt}+\int_{\Sigma }\alpha \widetilde{v}\psi
{\hskip1pt ds\hskip1pt dt}, \label{40-13}
\end{align}
\begin{align}
&-\int_{Q}\widetilde{\Phi }^{\lambda }(\psi _{1})_{t}{\hskip1pt dx\hskip1pt dt}+\int_{Q}\nabla
\widetilde{\Phi }^{\lambda }\cdot \nabla \psi _{1}{\hskip1pt dx\hskip1pt dt} 
\notag \\
&+\int_{Q}((\varphi _{\varepsilon ,\sigma }^{\lambda })^{2}+\varphi
_{\varepsilon ,\sigma }^{\lambda }\varphi _{\varepsilon ,\sigma }^{\ast
}+({\varphi}_{\varepsilon ,\sigma }^{\ast })^{2}-1)\widetilde{\Phi }^{\lambda
}\psi _{1}{\hskip1pt dx\hskip1pt dt}
=\int_{Q}\frac{\widetilde{\Theta }^{\lambda }}{\theta
_{\varepsilon ,\sigma }^{\lambda }{\theta}^{\ast }_{\varepsilon ,\sigma}}
\psi _{1}{\hskip1pt dx\hskip1pt dt}.  \label{40-14}
\end{align}
Now, {we observe that}
\begin{equation*}
\frac{\beta (\theta _{\varepsilon ,\sigma }^{\lambda })-\beta (\theta
_{\varepsilon ,\sigma }^{\ast })}{\lambda }=\beta ^{\prime }(\,\overline{\theta}_{\lambda } )\widetilde{\Theta }^{\lambda },
\end{equation*}
where $\overline{\theta}_{\lambda }$ {is a measurable function taking 
intermediate values between $\theta _{\varepsilon ,\sigma
}^{\ast }$ and $\theta _{\varepsilon ,\sigma }^{\lambda }$, a.e. in $Q$.
Moreover, due to \eqref{pier15} we have that} $\overline{\theta}_{\lambda }\rightarrow
{\theta}^{\ast }_{\varepsilon ,\sigma}$ strongly in $L^{2}(Q)$ as $\lambda
\rightarrow 0.$ Therefore, {by the Lipschitz continuity of some restriction of $\beta$ to a bounded subset of $(0,+\infty)$} we deduce that
\begin{equation*}
\beta ^{\prime }(\overline{\theta}_{\lambda })\rightarrow \beta ^{\prime
}({\theta}^{\ast }_{\varepsilon ,\sigma})\ \mbox{ strongly in }L^{2}(Q), \, \mbox{
as }\lambda \rightarrow 0,
\end{equation*}
{whence
\begin{equation*}
\frac{\beta (\theta _{\varepsilon ,\sigma }^{\lambda })-\beta (\theta
_{\varepsilon ,\sigma }^{\ast })}{\lambda }\rightarrow \beta ^{\prime
}({\theta}^{\ast }_{\varepsilon ,\sigma})\widetilde{Y} \ \mbox{ weakly
in }L^{1}(Q),\, \mbox{ as }\lambda \rightarrow 0,
\end{equation*}
first, and then weakly in $L^{2}(Q)$ due to the boundedness of $\beta ^{\prime }(\,\overline{\theta}_{\lambda } )\widetilde{\Theta }^{\lambda }$ in $L^2 (Q) $. 
Analogously, in view of \eqref{pier15},we have that 
\begin{align*}
((\varphi _{\varepsilon ,\sigma }^{\lambda })^{2}+\varphi
_{\varepsilon ,\sigma }^{\lambda }\varphi _{\varepsilon ,\sigma }^{\ast
}+({\varphi}_{\varepsilon ,\sigma }^{\ast })^{2}-1)\widetilde{\Phi } \, \to \, 
(3 (\varphi _{\varepsilon ,\sigma }^{\ast})^2 -1)\widetilde{\Phi }  
\ \mbox{ and}\\
\frac{\widetilde{\Theta }^{\lambda }}{\theta
_{\varepsilon ,\sigma }^{\lambda }{\theta}^{\ast }_{\varepsilon ,\sigma}} 
\to \frac{\widetilde{Y}}{({\theta}^{\ast }_{\varepsilon ,\sigma})^2 } \ \mbox{ weakly
in }L^{2}(Q),\, \mbox{ as }\lambda \rightarrow 0.
\end{align*}
Now, we can} pass to the limit in (\ref{40-13})-(\ref{40-14}) and {find out that}
\begin{align*}
-\int_{Q}\widetilde{Y}\psi _{t}{\hskip1pt dx\hskip1pt dt}-\int_{Q}\beta ^{\prime }(\theta
_{\varepsilon ,\sigma }^{\ast })\widetilde{Y}\Delta \psi {\hskip1pt dx\hskip1pt dt}-\int_{Q}
\widetilde{\Phi }\psi _{t}{\hskip1pt dx\hskip1pt dt} 
=\int_{Q}\widetilde{u}\psi {\hskip1pt dx\hskip1pt dt}+\int_{\Sigma }\alpha \widetilde{v}\psi
{\hskip1pt ds\hskip1pt dt},
\end{align*}
\begin{align*}
-\int_{Q}\widetilde{\Phi }(\psi _{1})_{t}{\hskip1pt dx\hskip1pt dt}+\int_{Q}\nabla \widetilde{%
\Phi }\cdot \nabla \psi _{1}{\hskip1pt dx\hskip1pt dt}+\int_{Q}(3(\varphi _{\varepsilon ,\sigma
}^{\ast })^{2}-1)\widetilde{\Phi }\psi _{1}{\hskip1pt dx\hskip1pt dt}\\
=\int_{Q}\frac{\widetilde{Y}}{(\theta _{\varepsilon ,\sigma }^{\ast
})^{2}}\psi _{1}{\hskip1pt dx\hskip1pt dt},
\end{align*}
which means that $\widetilde{Y},$ $\widetilde{\Phi }$ {yield a 
solution to (\ref{35})-(\ref{40}) (see~\eqref{40-6}--\eqref{pier26}).  
Since this solution is unique we obtain
$\widetilde{Y}=Y,$ $\widetilde{\Phi }=\Phi $
and it is the whole family to converge in \eqref{pier21} as $\lambda \to 0$.}
\hfill $\square $

\bigskip

{Next, let} us denote by $p_{\varepsilon ,\sigma }$ and $q_{\varepsilon
,\sigma }$ the dual variables and introduce the dual system%
\begin{equation}
(p_{\varepsilon ,\sigma })_{t}+\beta ^{\prime }(\theta _{\varepsilon ,\sigma
}^{\ast })\Delta p_{\varepsilon ,\sigma }+\frac{1}{(\theta _{\varepsilon
,\sigma }^{\ast })^{2}}q_{\varepsilon ,\sigma }=-I_{1,\varepsilon }^{\sigma
}, \ \mbox{ in }Q,  \label{41}
\end{equation}
\begin{equation}
(q_{\varepsilon ,\sigma })_{t}+\Delta q_{\varepsilon ,\sigma }-(3(\varphi
_{\varepsilon ,\sigma }^{\ast })^{2}-1)q_{\varepsilon ,\sigma
}+(p_{\varepsilon ,\sigma })_{t}=-I_{2,\varepsilon }^{\sigma }, \ \mbox{ in }
Q,  \label{41-1}
\end{equation}
\begin{equation}
{\frac{\partial p_{\varepsilon ,\sigma }}{\partial \nu }+\alpha (x)p_{
\varepsilon ,\sigma }=0, \ \, \frac{\partial q_{\varepsilon ,\sigma
}}{\partial \nu }=0, \ \mbox{ on }\Sigma ,}  \label{41-3}
\end{equation}
\begin{equation}
{p_{\varepsilon ,\sigma }(T)=0, \ \, q_{\varepsilon ,\sigma }(T)=0, \
\mbox{ in }\Omega ,}  \label{42}
\end{equation}
where%
\begin{align*}
{I_{1,\varepsilon }^{\sigma } =2\lambda _{1}(\theta _{\varepsilon
,\sigma }^{\ast }-\theta_{f})+\frac{1}{\varepsilon }(\xi_{\varepsilon
,\sigma }^{\ast }-\eta _{\varepsilon ,\sigma }^{\ast }), \quad
I_{2,\varepsilon }^{\sigma } =2\lambda _{2}(\varphi _{\varepsilon ,\sigma
}^{\ast }-\eta _{\varepsilon ,\sigma }^{\ast }), \quad \xi_{\varepsilon
,\sigma }^{\ast } =j_{\sigma }^{\prime }(\theta _{\varepsilon ,\sigma
}^{\ast }).}
\end{align*}
We note that
\begin{equation*}
I_{1,\varepsilon }^{\sigma },\, I_{2,\varepsilon }^{\sigma
},\,\xi_{\varepsilon ,\sigma }^{\ast }\in L^{\infty }(Q).
\end{equation*}

\medskip

\noindent \textbf{Proposition 4.4. }\textit{Assume }(\ref{10-0}
), (\ref{PF9}), (\ref{17})\textit{\ and }(\ref{20})\textit{. Then }
\textit{the dual system} (\ref{41})--(\ref{42}) \textit{has a
unique solution
{$(p_{\varepsilon ,\sigma },
q_{\varepsilon ,\sigma })$ with}}
\begin{equation}
p_{\varepsilon ,\sigma },\,{q_{\varepsilon ,\sigma } \in
L^{2}(0,T;H^{2}(\Omega )) \cap L^{\infty}(0,T;V)\cap
W^{1,2}([0,T];L^{2}(\Omega ))}  \label{43-1}
\end{equation}
\textit{and {such that the estimates}}
\begin{align}
& \left\Vert (p_{\varepsilon ,\sigma })_{t}\right\Vert _{L^{2}(Q)}^{2}+%
{\left\Vert p_{\varepsilon ,\sigma }\right\Vert _{L^\infty(0,T;V)}^{2}}
+\int_{0}^{T}\left\Vert \Delta p_{\varepsilon ,\sigma }(t)\right\Vert
_{L^{2}(\Omega )}^{2}{dt}\leq C,  \label{43-3} \\
& \left\Vert (q_{\varepsilon ,\sigma })_{t}\right\Vert
_{L^{2}(Q)}^{2}+ {}{\left\Vert q_{\varepsilon ,\sigma }\right\Vert_{
{L^\infty(0,T;V)}}^{2}} +\int_{0}^{T}\left\Vert \Delta
q_{\varepsilon ,\sigma }(t)\right\Vert _{L^{2}(\Omega )}^{2}{dt}\leq C
\label{pier7}
\end{align}
\textit{{hold \textit{independently of }$\sigma >0$}.}

\medskip

\noindent \textbf{Proof. }First, we make in (\ref{41})--(\ref{42}) the
variable transformation $t^{\prime }=T-t.$
{Then, the thesis follows
from} Proposition~4.2, by setting in the transformed system%
{
\begin{align*}
&a =\beta ^{\prime }({\theta}^{\ast }_{\varepsilon ,\sigma}),
\quad  b=-\frac{1}{({\theta}^{\ast }_{\varepsilon ,\sigma})^{2}},
\quad\omega =I_{1,\varepsilon }^{\sigma }, \\
&c =(3({\varphi}_{\varepsilon ,\sigma }^{\ast })^{2}-1),\quad
d =-1,\quad g = I_{2,\varepsilon }^{\sigma }, \ \hbox{ and }\,  \gamma =0.
\end{align*}
\hfill The estimates \eqref{43-3}--\eqref{pier7} can be obtained by standard computations, testing \eqref{41} by $- (p_{\varepsilon ,\sigma })_{t}$ and
\eqref{41-1} by $- (q_{\varepsilon ,\sigma })_{t}$, integrating, combining the resulting equalities, and so on. Finally, a comparison of terms in \eqref{41}
and \eqref{41-1} yields the desired estimates also for $ \Delta p_{\varepsilon ,\sigma } $ and $ \Delta q_{\varepsilon ,\sigma } $}
. \hfill $\square $

\medskip

\noindent \textbf{Proposition 4.5. }\textit{Under the
assumptions }(\ref{10-0}), (\ref{PF9}), (\ref{17})\textit{\ and }(\ref{20}),%
\textit{\ the optimality conditions for }$(P_{\varepsilon ,\sigma
})$ \textit{read}{
\begin{align}
&-\left( p_{\varepsilon ,\sigma }+2(u_{\varepsilon ,\sigma }^{\ast
}-u_{\varepsilon }^{\ast })\right)  \in \partial I_{K_{1}}(u_{\varepsilon
,\sigma }^{\ast }), \label{pier31} \\
&-\left( \alpha p_{\varepsilon ,\sigma }+2(v_{\varepsilon ,\sigma }^{\ast
}-v_{\varepsilon }^{\ast })\right)  \in \partial I_{K_{2}}(v_{\varepsilon
,\sigma }^{\ast }), \label{43-2}\\
&-\left( I_{3,\varepsilon }^{\sigma }+2(\eta _{\varepsilon ,
\sigma }^{\ast
}-\eta _{\varepsilon }^{\ast })\right)  \in \partial I_{K_{[-1,1]}}(\eta
_{\varepsilon ,\sigma }^{\ast }),   \label{pier32}
\end{align}
}
\textit{where}
\begin{equation}
{I_{3,\varepsilon }^{\sigma }=-2\lambda _{2}(\varphi _{\varepsilon ,\sigma
}^{\ast }-\eta _{\varepsilon ,\sigma }^{\ast })+\frac{1}{\varepsilon }
(\theta _{c}-{\theta}^{\ast }_{\varepsilon ,\sigma}).}  \label{46}
\end{equation}

\medskip

\noindent \textbf{Proof. }{Due to \eqref{41}--\eqref{42}, it is straightforward to realize that one can take $\psi =  p_{\varepsilon ,\sigma }$ in \eqref{40-6} 
(cf.~\eqref{40-7}) and $\psi_1 = q_{\varepsilon ,\sigma }$ in \eqref{pier26}. Then, by adding the equalities and integrating by parts in one term, we obtain}
\begin{align*}
&-\int_{Q}\left\{ (p_{\varepsilon ,\sigma })_{t}+\beta ^{\prime }(\theta
_{\varepsilon ,\sigma }^{\ast })\Delta p_{\varepsilon ,\sigma }+\frac{1}{
({\theta}^{\ast }_{\varepsilon ,\sigma})^{2}}q_{\varepsilon ,\sigma
}\right\} Y\,{dx\hskip1pt dt} \\
&-\int_{Q}\left\{ (q_{\varepsilon ,\sigma })_{t}+\Delta q_{\varepsilon
,\sigma }-(3({\varphi}_{\varepsilon ,\sigma }^{\ast })^{2}-1)q_{\varepsilon
,\sigma }+(p_{\varepsilon ,\sigma })_{t}\right\} \Phi \,{dx\hskip1pt dt}
\\
&=\int_{Q}\widetilde{u}\,p_{\varepsilon ,\sigma }\,{dx\hskip1pt dt} 
+ \int_{\Sigma}\alpha \,\widetilde{v}\,p_{\varepsilon ,\sigma }\,{ds\hskip1pt dt} .
\end{align*}
{Hence, with the help of (\ref{41})--(\ref{41-1}) we have}
\begin{equation}
\int_{Q}I_{1,\varepsilon }^{\sigma }Y\,{dx\hskip1pt dt}
+\int_{Q}I_{2,\varepsilon }^{\sigma }\Phi \,{dx\hskip1pt dt} 
= \int_{Q}\widetilde{u}\,p_{\varepsilon ,\sigma }\,{dx\hskip1pt dt} 
+ \int_{\Sigma}\alpha \, \widetilde{v}\,p_{\varepsilon ,\sigma }\,{ds\hskip1pt dt} .
\label{44}
\end{equation}
Then we write the optimality condition
\begin{equation*}
J_{\varepsilon ,\sigma }(u_{\varepsilon ,\sigma }^{\ast },v_{\varepsilon
,\sigma }^{\ast },\eta _{\varepsilon ,\sigma }^{\ast })\leq J_{\varepsilon
,\sigma }(\widehat{u},\widehat{v},\widehat{\eta }),\ \mbox{ for any }(%
\widehat{u},\widehat{v},\widehat{\eta })\in K_{1}\times K_{2}\times
K_{[-1,1]}.
\end{equation*}
In particular, {taking} $\widehat{u}=u_{\varepsilon ,\sigma }^{\lambda
},$ $\widehat{v}=v_{\varepsilon ,\sigma }^{\lambda },$ $\widehat{\eta }=\eta
_{\varepsilon ,\sigma }^{\lambda },$ making some computations, dividing by $
\lambda $ and letting $\lambda $ go to $0$ {lead to the inequality}
\begin{align*}
&2\lambda _{1}\int_{Q}({\theta}^{\ast }_{\varepsilon ,\sigma}-\theta _{f})Y
{\hskip1pt dx\hskip1pt dt} +2\lambda _{2}\int_{Q}(\varphi _{\varepsilon ,\sigma
}^{\ast }-\eta _{\varepsilon ,\sigma }^{\ast })(\Phi -\widetilde{\eta })
{\hskip1pt dx\hskip1pt dt} \\
&+\frac{1}{\varepsilon }\int_{Q}(\xi_{\varepsilon ,\sigma }^{\ast }Y+%
\widetilde{\eta }\theta _{c}-\widetilde{\eta }\theta _{\varepsilon ,\sigma
}^{\ast }-\eta _{\varepsilon ,\sigma }^{\ast }Y){\hskip1pt dx\hskip1pt dt} \\
&+2\int_{Q}(u_{\varepsilon ,\sigma }^{\ast }-u_{\varepsilon }^{\ast })%
\widetilde{u}\,{dx\hskip1pt dt} 
+{2} \int_{\Sigma }(v_{\varepsilon
,\sigma }^{\ast }-v_{\varepsilon }^{\ast })\widetilde{v}\,{ds\hskip1pt
dt} +2\int_{Q}(\eta _{\varepsilon ,\sigma }^{\ast }-\eta _{\varepsilon
}^{\ast })\widetilde{\eta }\,{dx\hskip1pt dt} \geq 0.
\end{align*}
With the previous {notation,} this yields
\begin{align}
&\int_{Q}I_{1,\varepsilon }^{\sigma }Y\,{dx\hskip1pt dt}
+\int_{Q}I_{2,\varepsilon }^{\sigma }\Phi \,{dx\hskip1pt dt} +\int_{Q} %
{I_{3,\varepsilon }^{\sigma }\widetilde{\eta }}\,{dx\hskip1pt dt}
\notag \\
&+2\int_{Q}(u_{\varepsilon ,\sigma }^{\ast }-u_{\varepsilon }^{\ast })%
\widetilde{u}\,{dx\hskip1pt dt} + {2} \int_{\Sigma
}(v_{\varepsilon ,\sigma }^{\ast }-v_{\varepsilon }^{\ast })\widetilde{v}\,%
{\hskip1pt ds\hskip1pt dt} +2\int_{Q}(\eta _{\varepsilon ,\sigma }^{\ast }-\eta
_{\varepsilon }^{\ast })\widetilde{\eta }\,{dx\hskip1pt dt} \geq 0.
\label{45}
\end{align}
{By comparing (\ref{44}) and (\ref{45}) we easily obtain}
\begin{align*}
\int_{Q}\widetilde{u}\left\{ p_{\varepsilon ,\sigma }+2(u_{\varepsilon
,\sigma }^{\ast }-u_{\varepsilon }^{\ast })\right\} {dx\hskip1pt dt}
&+\int_{\Sigma }\widetilde{v}\left\{ \alpha p_{\varepsilon ,\sigma
}+2(v_{\varepsilon ,\sigma }^{\ast }-v_{\varepsilon }^{\ast })\right\} 
{ ds\hskip1pt dt} \\
&+\int_{Q}\widetilde{\eta }\left\{ I_{3,\varepsilon }^{\sigma }+2(\eta
_{\varepsilon ,\sigma }^{\ast }-\eta _{\varepsilon }^{\ast })\right\}
{dx\hskip1pt dt} \geq 0.
\end{align*}
{Therefore, r}ecalling (\ref{34}) we finally have
\begin{align*}
&\int_{Q}(u_{\varepsilon ,\sigma }^{\ast }-u)\left\{ -\left( p_{\varepsilon
,\sigma }+2(u_{\varepsilon ,\sigma }^{\ast }-u_{\varepsilon }^{\ast
})\right) \right\} {dx\hskip1pt dt} \\
&+\int_{\Sigma }(v_{\varepsilon ,\sigma }^{\ast }-v)\left\{ -\left( \alpha
p_{\varepsilon ,\sigma }+2(v_{\varepsilon ,\sigma }^{\ast }-v_{\varepsilon
}^{\ast })\right) \right\} {ds\hskip1pt dt} \\
&+\int_{Q}(\eta _{\varepsilon ,\sigma }^{\ast }-\eta )\left\{ -\left(
I_{3,\varepsilon }^{\sigma }+2(\eta _{\varepsilon ,\sigma }^{\ast }-\eta
_{\varepsilon }^{\ast })\right) \right\}{dx\hskip1pt dt} \geq 0,
\end{align*}
for any $(u,v,\eta )\in K_{1}\times K_{2}\times K_{[-1,1]}.$ This implies (%
\ref{43-2}), as claimed. \hfill $\square $

\medskip

\noindent \textbf{Theorem 4.6. }\textit{ Assume}
(\ref{10-0}), (\ref{PF9}), (\ref{17}), (\ref{20}) \textit{and let} 
$\{(u_{\varepsilon }^{\ast },v_{\varepsilon }^{\ast },\eta
_{\varepsilon }^{\ast }),(\theta _{\varepsilon }^{\ast },\varphi
_{\varepsilon }^{\ast })\}$ \textit{be optimal in }$(P_{\varepsilon }).$
\textit{Then, the optimality conditions for }$(P_{\varepsilon })$ \textit{%
read\ }{
\begin{equation}
\left\{
\begin{array}{cl}
u_{\varepsilon }^{\ast }(t,x)= u_{m}, &\mbox{ \textit{on} \,}\{(t,x)\in Q:\
p_{\varepsilon }(t,x)>0\} \\
u_{m} \leq u_{\varepsilon }^{\ast }(t,x)\leq u_{M}, 
&\mbox{ \textit{on} \,}\{(t,x)\in Q:\ p_{\varepsilon }(t,x)=0\} \\
u_{\varepsilon }^{\ast }(t,x) = u_{M}, &\mbox{ \textit{on} \,}\{(t,x)\in Q:\
p_{\varepsilon }(t,x)<0\}
\end{array}
\right. \ ,  \label{u*opt}
\end{equation}
\begin{equation}
\left\{
\begin{array}{cl}
v_{\varepsilon }^{\ast }(t,x)= v_{m}, &\mbox{ \textit{on} \,}\{(t,x)\in \Sigma:\
p_{\varepsilon }(t,x)>0\} \\
v_{m} \leq v_{\varepsilon }^{\ast }(t,x)\leq v_{M}, 
&\mbox{ \textit{on} \,}\{(t,x)\in \Sigma:\ p_{\varepsilon }(t,x)=0\} \\
v_{\varepsilon }^{\ast }(t,x) = v_{M}, &\mbox{ \textit{on} \,}\{(t,x)\in \Sigma:\
p_{\varepsilon }(t,x)<0\}
\end{array}
\right. \ ,  \label{v*opt}
\end{equation}
\begin{equation}
\left\{
\begin{array}{cl}
\eta_{\varepsilon }^{\ast }(t,x)= -1, &\mbox{ \textit{on} \,}\{(t,x)\in Q:\
 I_{3,\varepsilon } (t,x)>0\} \\
-1 \leq \eta_{\varepsilon }^{\ast }(t,x)\leq 1, 
&\mbox{ \textit{on} \,}\{(t,x)\in Q:\  I_{3,\varepsilon } (t,x)=0\} \\
\eta_{\varepsilon }^{\ast }(t,x) = 1, &\mbox{ \textit{on} \,}\{(t,x)\in Q:\
 I_{3,\varepsilon } (t,x)<0\}
\end{array}
\right. \ ,  \label{eta*opt}
\end{equation}}

\smallskip
\noindent
\textit{where} {$(p_{\varepsilon }, q_{\varepsilon }) $} \textit{is the
solution to {the problem}}
\begin{equation}
(p_{\varepsilon })_{t}+\beta ^{\prime }(\theta _{\varepsilon }^{\ast
})\Delta p_{\varepsilon }+\frac{1}{(\theta _{\varepsilon }^{\ast })^{2}}
q_{\varepsilon }=-I_{1,\varepsilon },\ \mbox{ \textit{in} }Q,  \label{d1}
\end{equation}
\begin{equation}
(q_{\varepsilon })_{t}+\Delta q_{\varepsilon }-(3(\varphi _{\varepsilon
,\sigma }^{\ast })^{2}-1)q_{\varepsilon }+(p_{\varepsilon
})_{t}=-I_{2,\varepsilon },\ \mbox{ \textit{in} }Q,  \label{d2}
\end{equation}
\begin{equation}
\frac{\partial p_{\varepsilon }}{\partial \nu }+\alpha (x)p_{\varepsilon
}=0,\ \ \frac{\partial q_{\varepsilon }}{\partial \nu }=0,\
\mbox{ \textit{on} }\Sigma ,
\end{equation}
\begin{equation}
p_{\varepsilon }(T)=0,\ \ q_{\varepsilon }(T)=0,\ \mbox{ \textit{in} }\Omega
,  \label{d3}
\end{equation}
\textit{and {where}}
\begin{align*}
& {I_{1,\varepsilon }}=2\lambda _{1}(\theta _{\varepsilon }^{\ast
}-\theta _{f})+\frac{1}{\varepsilon }(\xi _{\varepsilon }^{\ast }-\eta
_{\varepsilon }^{\ast }),\ \ \mbox{\textit{{with }}}\ \xi _{\varepsilon
}^{\ast }\in \partial j(\theta _{\varepsilon }^{\ast })\ \,%
\mbox{\textit{{a.e. in }}}Q, \\
& I_{2,\varepsilon }=2\lambda _{2}(\varphi _{\varepsilon }^{\ast }-\eta
_{\varepsilon }^{\ast }),\quad \ I_{3,\varepsilon }=-2\lambda _{2}(\varphi
_{\varepsilon }^{\ast }-\eta _{\varepsilon }^{\ast })+\frac{1}{\varepsilon }
(\theta _{c}-\theta _{\varepsilon }^{\ast }).
\end{align*}

\medskip

\noindent \textbf{Proof. }Under the hypotheses, problem $(P_{\varepsilon
,\sigma })$ has a
{minimizer $(u_{\varepsilon ,\sigma }^{\ast },v_{\varepsilon
,\sigma }^{\ast },\eta _{\varepsilon ,\sigma }^{\ast })$\textit{\ }with the
corresponding pair \ (${\theta}^{\ast }_{\varepsilon ,\sigma},\varphi
_{\varepsilon ,\sigma }^{\ast })$ solving \eqref{PF1}--\eqref{PF6}}. We pass
to the limit in $(P_{\varepsilon ,\sigma }).$ According to Proposition~4.1%
{, we have the convergences in  {\eqref{pier2}--\eqref{pier3} and 
\eqref{pier8}--\eqref{pier9}, in which however the actual limits are 
$\theta _{\varepsilon }^{\ast }$ and $\varphi_{\varepsilon }^{\ast }$}. 
By the estimates (\ref{43-3})--(\ref{pier7}) in Proposition~4.4, 
at least for a subsequence we have~that}
\begin{align*}
p_{\varepsilon ,\sigma }\rightarrow p_{\varepsilon}, \, \ q_{\varepsilon
,\sigma }\rightarrow q_{\varepsilon} \mbox{
weakly in }L^{2}(0,T;H^{2}(\Omega ))\cap W^{1,2}([0,T];L^{2}(\Omega )),
\qquad\qquad  \notag \\
\mbox{{weakly*} in } L^{\infty }(0,T;V),\mbox{ and strongly in }
L^{2}(0,T;V),\mbox{ \textit{as} }{\sigma \rightarrow 0} .
\end{align*}
Then, {recalling \eqref{pier31} and passing to the limit we find that
\begin{equation*}
-\left( p_{\varepsilon ,\sigma
}+2(u_{\varepsilon ,\sigma }^{\ast }-u_{\varepsilon }^{\ast })\right)
\rightarrow -p_{\varepsilon }\ \mbox{ strongly in }L^{2}(Q),\mbox{ as }
\sigma \rightarrow 0,
\end{equation*}
which, along with \eqref{pier2}, yields
\begin{equation}
-p_{\varepsilon }\in \partial I_{K_{1}}(u_{\varepsilon }^{\ast }) {,
\quad{\hbox{a.e. in } Q }}  \label{48}
\end{equation}
for $ \partial I_{K_{1}}$ is maximal monotone and so strongly-weakly closed.} 
The same argument works for the other two controllers {in \eqref{43-2} 
and \eqref{pier32},} hence we obtain
\begin{equation}
-\alpha p_{\varepsilon }\in \partial I_{K_{2}}(v_{\varepsilon }^{\ast }) %
{, \quad{\hbox{a.e. on } \Sigma, }}  \label{49}
\end{equation}
and
\begin{equation}
- I_{3,\varepsilon }\in \partial I_{K_{[-1,1]}} (
\eta _{\varepsilon }^{\ast }) %
{, \quad{\hbox{a.e. in } Q , }}  \label{50}
\end{equation}
because {of the convergence}
\begin{multline*}
{-( I_{3,\varepsilon }^{\sigma } + 2(\eta _{\varepsilon
,\sigma }^{\ast }-\eta _{\varepsilon }^{\ast }))
=2\lambda _{2}(\varphi _{\varepsilon ,\sigma
}^{\ast }-\eta _{\varepsilon ,\sigma }^{\ast })-\frac{1}{\varepsilon }
(\theta _{c}-{\theta}^{\ast }_{\varepsilon ,\sigma})-2(\eta _{\varepsilon
,\sigma }^{\ast }-\eta _{\varepsilon }^{\ast })}\\
{\rightarrow \ 2\lambda _{2}(\varphi _{\varepsilon }^{\ast }-\eta
_{\varepsilon }^{\ast })- \frac{1}{\varepsilon }(\theta _{c}-\theta
_{\varepsilon }^{\ast })=-I_{3,\varepsilon }\ \mbox{
strongly in }L^{2}(Q), \mbox{ as }\sigma \rightarrow 0.}
\end{multline*}
It is easily seen that
\begin{multline*}
-I_{2,\varepsilon }^{\sigma }=2\lambda _{2}(\varphi _{\varepsilon ,\sigma
}^{\ast }-\eta _{\varepsilon ,\sigma }^{\ast }) \rightarrow -2\lambda
_{2}(\varphi _{\varepsilon }^{\ast }-\eta _{\varepsilon }^{\ast
})=-I_{2,\varepsilon }\ \mbox{ strongly in }L^{2}(Q), \mbox{ as }\sigma
\rightarrow 0.
\end{multline*}
{Letting $\xi_{\varepsilon }^{\ast }$ denote the weak* limit in $L^\infty (Q)$ 
of some subsequence of $\{\xi_{\varepsilon ,\sigma }^{\ast }\}$, it turns out that
$\xi_{\varepsilon }^{\ast }\in \partial j(\theta _{\varepsilon }^{\ast
})$ {a.e. in} $Q.$ Indeed, recalling that $\xi_{\varepsilon ,\sigma
}^{\ast }=j_{\sigma }^{\prime }({\theta}^{\ast }_{\varepsilon ,\sigma}),$ we
can write
\begin{equation*}
\int_{Q}(j_{\sigma }({\theta}^{\ast }_{\varepsilon ,\sigma})-j_{\sigma
}(z))\,{dx\hskip1pt dt} \leq \int_{Q}j_{\sigma }^{\prime }(\theta
_{\varepsilon ,\sigma }^{\ast })({\theta}^{\ast }_{\varepsilon ,\sigma}-z)\,%
{\hskip1pt dx\hskip1pt dt}
\end{equation*}
for any $z\in L^{2}(Q),$ and pass to the limit as $\sigma
\rightarrow 0$ taking {(\ref{33}) into account. Thus, we deduce that}
\begin{equation*}
\int_{Q}(j(\theta _{\varepsilon }^{\ast })-j(z))\,{dx\hskip1pt dt} \leq
\int_{Q}\xi _{\varepsilon }^{\ast }(\theta _{\varepsilon }^{\ast }-z)\,%
{\hskip1pt dx\hskip1pt dt} ,
\end{equation*}
which implies $\xi_{\varepsilon }^{\ast }\in \partial j(\theta _{\varepsilon
}^{\ast })$ {a.e. in} $Q.$ Consequently, we have that
\begin{multline*}
I_{1,\varepsilon }^{\sigma }=2\lambda _{1}(\theta _{\varepsilon ,\sigma
}^{\ast }-\theta_{f})+\frac{1}{\varepsilon }(\xi_{\varepsilon ,\sigma
}^{\ast }-\eta _{\varepsilon ,\sigma }^{\ast }) \\
\rightarrow \ 2\lambda _{1}(\theta _{\varepsilon }^{\ast }-\theta_{f})+\frac{%
1}{\varepsilon }(\xi _{\varepsilon }^{\ast }-\eta _{\varepsilon }^{\ast
})=I_{1,\varepsilon },\ \mbox{ weakly in }L^{2}(Q), \mbox{ as }\sigma
\rightarrow 0.
\end{multline*}
The above arguments prove} that the solution
to (\ref{41})--(\ref{42}) converges to the solution to 
(\ref{d1})--(\ref{d3}) as $\sigma \rightarrow 0$. 
{In fact, due to the uniform boundedness properties ensured by  
assumptions (\ref{17}) and (\ref{20}), we also point out that}
\begin{align*}
{\beta ^{\prime }(\theta _{\varepsilon ,\sigma }^{\ast }) \to
\beta ^{\prime }(\theta _{\varepsilon}^{\ast }) , \quad 
\frac{1}{(\theta _{\varepsilon ,\sigma }^{\ast })^{2}} \to 
\frac{1}{(\theta _{\varepsilon }^{\ast })^{2}} , \quad
(3(\varphi_{\varepsilon ,\sigma }^{\ast })^{2}-1) \to
(3(\varphi_{\varepsilon}^{\ast })^{2}-1)}\\
{\ \hbox{ weakly* in } L^\infty (Q) \hbox{ and strongly in } L^2 (Q), \hbox{ as }
\sigma \rightarrow 0.}
\end{align*}

{Note that the selection $\xi _{\varepsilon }^{\ast }$ from $\partial j(\theta _{\varepsilon}^{\ast })= H (\theta_{\varepsilon }^{\ast })$ which is present in $I_{1,\varepsilon} $ is not uniquely determined unless $\theta_{\varepsilon }^{\ast } \not= \theta_{c}$ a.e. in $Q$. On the other hand, the pair 
$(p_\varepsilon , q_\varepsilon) $ turns out to be the unique solution of the problem 
(\ref{d1})--(\ref{d3}) once $\xi _{\varepsilon }^{\ast }$ is fixed in $I_{1,\varepsilon}$.}

{Now, in order to conclude the proof it suffices to notice that, e.g., $\partial I_{K_{1}}(u_{\varepsilon }^{\ast
})$ is exactly $N_{K_{1}}(u_{\varepsilon }^{\ast })$, the normal cone to $K_{1}$
at $u_{\varepsilon }^{\ast }.$ Then, it is straightforward to derive 
(\ref{u*opt})--(\ref{eta*opt}), as claimed.} 
\hfill $\square $
\bigskip

\section*{Acknowledgements}
This research activity has been performed in the framework of an
Italian-Romanian collaboration leading to the approval of a three-year grant for the project on ``Nonlinear partial differential equations (PDE) with applications in modeling cell growth, chemotaxis and phase transition''. {Moreover, 
the financial support of the FP7-IDEAS-ERC-StG \#256872
(EntroPhase) is gratefully acknowledged by the authors. The present paper 
also benefits from the support of the MIUR-PRIN Grant 2010A2TFX2 ``Calculus of Variations'' for PC, the GNAMPA (Gruppo Nazionale per l'Analisi Matematica, la Probabilit\`a e le loro Applicazioni) of INdAM (Istituto Nazionale di Alta Matematica) for PC and ER, and the project CNCS-UEFISCDI, PN-II-ID-PCE-2011-3-
0027 for GM.}

\bigskip

\end{document}